\documentclass[12pt,reqno]{amsart}
\usepackage{amssymb,delarray}
\usepackage{amsfonts}
\usepackage{epsfig}
\usepackage[all]{xy}
\usepackage[utf8]{inputenc}
\usepackage{pstricks-add}

\def\leq{\leqslant}
\def\geq{\geqslant}

\def\l{\operatorname{l}}
\def\pin{{\operatorname{
cdt}}}

\newtheorem{thm}{Theorem}
\newtheorem{lem}
{Lemma}
\newtheorem{prop}
{Proposition}
\newtheorem{claim}
{Claim}

\newtheorem{cor}
{Corollary}
\newtheorem{rem}
{Remark}
{Question}

{\catcode`\@=11
\gdef\n@te#1#2{\leavevmode\vadjust{%
 {\setbox\z@\hbox to\z@{\strut#1}%
  \setbox\z@\hbox{\raise\dp\strutbox\box\z@}\ht\z@=\z@\dp\z@=\z@%
  #2\box\z@}}}
\gdef\leftnote#1{\n@te{\hss#1\quad}{}}
\gdef\rightnote#1{\n@te{\quad\kern-\leftskip#1\hss}{\moveright\hsize}}
\gdef\?{\FN@\qumark}
\gdef\qumark{\ifx\next"\DN@"##1"{\leftnote{\rm##1}}\else
 \DN@{\leftnote{\rm??}}\fi{\rm??}\next@}}

\begin{document}
\baselineskip=13.7pt plus 2pt 

\title[Covering semigroups] {Covering semigroups}
\author[Vik.S. Kulikov and V.M. Kharlamov]{Vik.S. Kulikov and V.M. Kharlamov}

\address{Steklov Mathematical Institute}
 \email{kulikov@mi.ras.ru}

 \address{
UFR de Math\'ematiques et IRMA \\ Universit\'e de Strasbourg et CNRS
}
\email{ kharlam@math.unistra.fr}
\dedicatory{} \subjclass{}
\thanks{The first author was partially supported by grants of NSh-4713.2010.1, RFBR
11-01-00185, and by AG Laboratory HSE, RF government grant, ag.
11.G34.31.0023. The second author was partially funded by the
ANR-09-BLAN-0039-01 grant of Agence Nationale de la Recherche, and
is a member of FRG Collaborative Research "Mirror Symmetry $\&$
Tropical Geometry" (Award No. 0854989).}


\begin{abstract}
We introduce and study a semigroup structure on the set of
irreducible components of the Hurwitz space of marked coverings of a
complex projective curve with given Galois group of the coverings
and fixed ramification type. As application, we give new conditions
on the ramification type that are sufficient for irreducibility of
the Hurwitz spaces, suggest some bounds on the number of
irreducibility components under certain more general conditions, and
show that the number of irreducible components coincides with the
number of topological classes of the coverings if the number of
brunch points is big enough.
\end{abstract}

\maketitle
\setcounter{tocdepth}{1}


\def\st{{\sf st}}



\section*{Introduction} \label{introduc}
Let $f:E\to F$ be a finite morphism between complex non-singular
irreducible projective curves. Denote by $\mathbb C(E)$ and
$\mathbb C(F)$ the fields of rational functions
on $E$ and $F$, respectively. The morphism $f$ defines a finite
extension $f^*: \mathbb C(F)\hookrightarrow \mathbb C(E)$ of the
field $\mathbb C(F)$ (reciprocally, the field extension defines the
covering $f$ uniquely up to isomorphisms of coverings over a fixed
base). We denote by $G$ the Galois group of the Galois closure of
this extension and call it  the {\it Galois group of $f$}.

Let us fix a point $q\in F$ that is not a branch point of $f$ and
order the points of $E$ lying over $q$. We call the morphism $f$ with a fixed ordering
of the points of $f^{-1}(q)$ a {\it marked covering}.

Consider the fundamental group $\pi_1(F\setminus B, q)$
of the complement of the branch set $B \subset F$ of a marked
covering $f$ of degree $d=\deg f$. Then, the ordering of the points of $f^{-1}(q)$
defines a homomorphism $f_*:\pi_1(F\setminus B, q)\to \mathcal
S_d$ of $\pi_1(F\setminus B, q)$ to the symmetric group $\mathcal
S_d$. Due to irreducibility of $E$, the image $\text{im}
f_*\subset\mathcal S_d$ acts transitively on $f^{-1}(q)$ and is
isomorphic to $G$, so that we can identify $\text{im} f_*$ and $G$
and thus fix this embedding $G\hookrightarrow\mathcal S_d$.

The movement along a standard simple loops $\gamma$ around branch
points $b\in B$
the {\it local monodromy} $f_*(\gamma)\in G$ of $f$ at $b$. The
homotopy class of this standard loop, and hence the local monodromy,
are defined by $b$ uniquely only up to conjugation, in
$\pi_1(F\setminus B, q)$ and $G$, respectively. We denote by $O
\subset G$ the union of the conjugacy classes of all the local
monodromies of $f$ and call the pair $(G,O)$ the {\it equipped
Galois group} associated with $f$. The collection
$\tau=(\tau_1C_1,\dots,\tau_mC_m)$, where $C_1,\dots, C_m$ list all
the conjugacy classes included in $O$ and $\tau_i$ counts the number
of branch points of $f$ with the local monodromies belonging to
$C_i$, is called the {\it monodromy type} of $f$.

The degree $d$  marked coverings of $F$ with Galois group $G$ and
monodromy type $\tau$ form a so called Hurwitz space
$\text{HUR}_{d,G,\tau} (F)$ (for precise definitions see subsection
\ref{Hursp}).

In the case $F=\mathbb P^1$, $G=\mathcal S_d$ and $O$ is the set of
transpositions, the famous Clebsch -- Hurwitz Theorem \cite{Cl},
\cite{H} states that  $\text{HUR}_{d,\mathcal S_d,\tau} (\mathbb
P^1)$ consists of a single irreducible component if $\tau= (nO)$
with even $n\geq 2(d-1)$ and it is empty otherwise. Generalizations
of Clebsch -- Hurwitz Theorem were obtained in \cite{BE}, \cite{W},
\cite{P}, \cite{FV}, and  \cite{Ku1} --  \cite{Ku3}. In particular,  Clebsch --
Hurwitz Theorem was extended to the following cases: in \cite{BE},
if all but one local monodromies are transpositions; in \cite{W}, if
all but two local monodromies are transpositions; in \cite{P}, if
all local monodromies are either transpositions or cyclic
permutations of length three; and in \cite{Ku1}, if there are at
least $3(d-1)$ transpositions among the local monodromies.

In \cite{Ku2}, it is proved that for an equipped group $(\mathcal
S_d,O)$ such that the first conjugacy class $C_1$ of $O$ contains an
odd permutation leaving fixed at least two elements, the Hurwitz
space $\text{HUR}_{d,\mathcal S_d,\tau} (\mathbb P^1)$ is
irreducible if $\tau_1$ is big enough. On the other hand, the
example in \cite{W} shows that $\text{HUR}_{8,\mathcal S_8,\tau}
(\mathbb P^1)$ consists of at least two irreducible components if
$\tau=(1C_1,1C_2,1C_3)$, where $C_1$ is the conjugacy class of the
permutation $(1,2)(3,4,5)$, $C_2$ is the conjugacy class of
$(1,2,3)(4,5,6,7)$, and $C_3$ is the conjugacy class of
$(1,2,3,4,5,6,7)$. Articles  \cite{FV} and \cite{Ku3} are devoted to partial generalizations
of Clebsch -- Hurwitz Theorem to the case of arbitrary group $G$. In particular, in \cite{Ku3}, it was
proved that for a fixed equipped finite group $(G,O)$ the number of
irreducible components of $\text{HUR}_{d,G,\tau} (\mathbb P^1)$ (if
it is non-empty) does not depend on $\tau$ if all $\tau_i$ are big
enough.

For higher genus, the irreducibility of $\text{HUR}_{d,\mathcal
S_d,\tau} (F)$ is proved in \cite{GHS} under hypothesis that $n\ge
2d$ and all local monodromies are transpositions. After that, this
result was improved, first, in \cite{Ka} where the hypothesis $n\ge
2d$ was replaced by $n\ge 2d-2$, and next, in \cite{V}, where the
second hypothesis was replaced by assumption that all but one local
monodromies are transpositions. Let us notice that the
irreducibility of the quotient of $\text{HUR}_{d,\mathcal S_d,\tau}
(F)$ by the action of the mapping class group of $F$ (considered as
a real surface) was proved in \cite{BE} under a weaker hypothesis
$n>\frac{d}2$.

One of the aims of this article is to extend results of \cite{Ku1}
-- \cite{Ku3} from $F=\mathbb P^1$ to the case of $F$ of arbitrary
genus. The approach used there for counting the number of
irreducible components of $\text{HUR}_{d,G,\tau}(\mathbb P^1)$ is
based on a systematic work with semigroups over groups; in
particular, factorization semigroups $S(G,O)$ with factors belonging
to $O$ (cf., subsections \ref{semig} and  \ref{semigr} below) play
the crucial role in this study, especially since subsets of elements
of type $\tau$ of subsemigroup $S(G,O)_{\bf 1}^G\subset S(G,O)$ are
in a canonical bijection with the sets of irreducible components of
the Hurwitz space $\text{HUR}_{d,G, \tau}(\mathbb P^1)$.

In the present paper, to treat the coverings of projective curves
(or, the same, real surfaces) of arbitrary genus we generalize the
notion of factorization semigroups to that of semigroups of marked
coverings. One can consider different levels of the equivalence
relations of coverings and so we introduce, respectively, different
species of semigroups of marked coverings. The equivalence relation
of the level that is most appropriate to construction of Hurwitz
spaces  is based essentially on moving of branch points, while that
the level most appropriate to topological classification of
coverings (like in \cite{BE}, for example) includes, in addition,
the  action on the base of coverings by the whole mapping class
group.  In particular, considering the coverings up to moving of
branch points we introduce a semigroup $G\mathbb S_{d}(G,O)$ of
marked degree $d$ coverings with Galois group $G$ and set of local
monodromies $O\subset G$. If we consider the same coverings up to
the action of the modular group, then we obtain another semigroup,
which we denote by $GW\mathbb S_{d}(G,O)$. They are related by a
natural epimorphism $\Phi :G\mathbb S_{d}(G,O)\to GW\mathbb
S_{d}(G,O)$ of semigroups. Similar to genus $0$ case, certain
specific subsemigroups of these two semigroups are in a canonical
bijection with the set of irreducible components of the Hurwitz
space $\text{HUR}_{d,G, 
}(F)$ and, respectively, the set of topological classes of marked
degree $d$ ramified coverings of $F$ with Galois groups $G$.

By definition, the {\it monodromy type} of an element $s=(f:E\to F)$
belonging to one of these semigroups is the collection
$\tau(s)=(\tau_1C_1,\dots, \tau_mC_m)$ of local monodromies of $f$.
The  monodromy type behaves additively and gives a homomorphism from
semigroups of coverings to the semigroup  $\mathbb Z^m_{\geq 0}$.
Therefore, for any constant $T\in \mathbb N$, there appear well
defined subsemigroups
$$G\mathbb S_{d, T}(G,O)=\{ s\in G\mathbb S_{d}(G,O)\mid
\tau_i(s)\geq T\, \, \text{for}\, \, i=1,\dots,m \}$$
and
$$
GW\mathbb S_{d,T}(G,O)=\{ s\in GW\mathbb S_{d}(G,O)\mid \tau_i(s)\geq
T\, \, \text{for}\, \, i=1,\dots,m \}.
$$

The main results are as follows.

\begin{thm}\label{cor3} For any equipped finite group $(G,O)$ such that the elements of $O$ generate the group $G$,
there is a constant $T\in \mathbb N$ such that the restriction of
$\Phi$ to $G\mathbb S_{d,T}(G,O)$ is an isomorphism of
$G\mathbb S_{d,T}(G,O)$ and $GW\mathbb S_{d,T}(G,O)$.
\end{thm}

In \cite{Ku3}, there was defined an {\it ambiguity index} for each
equipped finite group $(G,O)$ (see subsection \ref{sub3}).
\begin{thm} \label{TH2} For each equipped finite group $(G,O)$, $O=C_1\sqcup\dots
\sqcup C_m$, such that the elements of $O$ generate the group $G$,
there is a constant $T$ such that for any projective irreducible
non-singular curve $F$ the number of irreducible
components of each non-empty Hurwitz space $\text{HUR}_{d,G,
\tau}(F)$ is equal to $a_{(G,O)}$ if $\tau_i\geq T$ for
all $i=1,\dots,m$.

If the elements of $O_k=C_1\sqcup\dots \sqcup C_k$ for some $k<m$
generate the group $G$, then there is a constant $T'$ such that the
number of irreducible components of $\text{HUR}_{d,G, \tau}^m( F)$
is less or equal to $a_{(G,O_k)}$ if $\tau_i\geq T'$ for all
$i=1,\dots, k$.
\end{thm}

\begin{thm} \label{TH1} Let $C$ be the conjugacy class of an odd
permutation $\sigma\in \mathcal S_d$ such that $\sigma$ leaves fixed
at least two elements. Then there is a constant $N_C$ such that  for
any projective irreducible non-singular curve $F$ the
Hurwitz space {\rm $\text{HUR}_{d,\mathcal S_d, \tau}( F)$}
is irreducible if $C$ enters in $\tau$ with a factor $\ge N_C$.
\end{thm}

The article consists of  two sections. Section \ref{secA} is devoted
to the algebraic part of the proof. In subsections \ref{semig}  --
\ref{sub3} we fix notation and recall necessary definitions and
results from \cite{Ku1} --  \cite{Ku3}. In subsection \ref{adm} we
introduce a notion of admissible subgroups of the automorphism
groups of free groups, which is necessary for the next subsection
where we define the algebraic coverings semigroups. The remaining
subsections of Section \ref{secA} contain the proofs of the
algebraic part of main results. Section \ref{secG} starts from two
preliminary subsections where we introduce such auxiliary notions
like monodromy encoding of ramified coverings and skeletons of
surfaces.  In \ref{geom-free}  -- \ref{weak-semigroups} we introduce
a series of geometric covering semigroups and prove comparison
statements between algebraic and geometric covering semigroups. In
the final subsections, we relate elements of the geometric coverings
semigroups with irreducible components of Hurwitz spaces and
complete the proofs of main theorems.

\section{Semigroups over groups} \label{secA}
\subsection{Definition of semigroups over groups} \label{semig}
Here, we recall basic definitions and some properties of semigroups
over groups with a special emphasis to factorization semigroups (for
more details, see \cite{Ku} -- \cite {Ku3}).

A collection $(S,G,\alpha,\rho) $, where $S$ is a semigroup, $G$ is
a group, and $\alpha :S\to G$, $\rho :G\to \text{Aut}(S)$ are
homomorphisms, is called {\it a semigroup $S$ over a group $G$} if
for all $s_1,s_2\in S$ we have
\begin{equation} \label{relx}
s_1\cdot s_2=\rho (\alpha (s_1))(s_2)\cdot s_1= s_2\cdot\lambda
(\alpha (s_2))(s_1),\end{equation} where $\lambda(g)=\rho (g^{-1})$.

Let $(S_1,G_1,\alpha_1,\rho_1)$ and $(S_2,G_2,$ $\alpha_2,\rho_2)$
be two semigroups over groups $G_1$ and $G_2$, respectively. A pair
$(h_1,h_2)$ of homomorphisms $h_1:S_1\to S_2$ and $h_2: G_1\to G_2$
is called {\it a homomorphism of semigroups over groups} if
\begin{itemize}
\item[{($i$)}] $h_2\circ \alpha_{1}=\alpha_{2}\circ h_1$,
\item[($ii$)] $\rho_{2}(h_2(g))(h_1(s))=h_1(\rho_{1}(g)
(s))$ for all $s\in S_1$ and all $g\in G_1$.
\end{itemize}
In particular, if $G_1=G_2=G$, then a homomorphism of semigroups
$\varphi : S_1\to S_2$ is said {\it to be defined over} $G$ if
$\alpha_1(s)=\alpha_2(\varphi(s))$ and
$\rho_2(g)(\varphi(s))=\varphi(\rho_1(g)
(s))$ for all $s\in S_1$ and $g\in G$.

\subsection{Factorization semigroups} \label{semigr}
One of the main examples of semigroups over groups is given by, so
called, factorization semigroups. To define them, consider an {\it
equipped group} $(G,O)$, that is,  $G$ is a group and $O$ is a
subset of $G$ invariant under the inner automorphisms. Here and
further on, we assume that:
\begin{enumerate}
\item[(i)] ${\bf{1}}\not\in O$;
\item[(ii)] $O$ consists of a
finite number of conjugacy classes $C_i$ of $G$, $O= C_1\sqcup
\dots\sqcup C_m$;
\item[(iii)] the (linear) ordering of these conjugacy classes is fixed.
\end{enumerate}
By {\it homomorphisms of equipped groups} $(G_1,O_1)$ and
$(G_2,O_2)$ we understand homomorphisms $f: G_1\to G_2$ such that
$f(O_1)\subset O_2$.

{\it The factorization semigroup with factors in} $O$  is, by
definition, the semigroup $S(G,O)$ generated by alphabet $\mathcal
X_O=\{ x_g \mid g\in O\}$ and subject to relations
\begin{equation} \label{rel1} x_{g_1}\cdot
x_{g_2}=x_{g_2}\cdot\,x_{g_1^{g_2}}, \qquad g_1,g_2\in O,
\end{equation} where $g_1^{g_2}$
denotes $g_2^{-1}g_1g_2$. The homomorphism $\alpha :S(G,O)\to G$,
given by $\alpha (x_g)=g$ for each $x_g\in \mathcal X_O$, is called
the {\it product homomorphism}. {\it The simultaneous conjugation}
$$
x_a\in \mathcal X_O\mapsto x_{gag^{-1}}\in \mathcal X_O
$$
defines a homomorphism $G\to \text{Aut} (S(G,O))$, which we denote
by $\rho$. It is easy to see that under such a choice,
$(S(G,O),G,\alpha,\rho )$ is a semigroup over $G$.

Note that there is a well defined {\it length homomorphism} of
semigroups,
$$
\l :S(G,O)\to \mathbb Z_{\geq 0}=\{ {{a}}\in \mathbb Z\mid {{a}}\geq
0\}
$$
that is defined by $\l(x_{g_1}\cdot .\, .\, .\, \cdot x_{g_n})=n$.

Put $\rho _{S} =\rho \circ\alpha, \lambda _{S} =\lambda
\circ\alpha$, where as above $\lambda (g)=\rho(g^{-1})$.

\begin{claim} {\rm (\cite{KK})} \label{cl1}
For all $s_1,\, s_2\in S(G,O)$ we have
\[
s_1\cdot s_2=s_2\cdot \lambda_S(s_2)(s_1)=\rho_S(s_1)(s_2)\cdot s_1.
\]
\end{claim}

To each $s=x_{g_1}\cdot\, .\, . \, . \, \cdot x_{g_n}\in S(G,O)$, we
associate a subgroup $G_s$ of $G$ generated by the images $\alpha
(x_{g_1})=g_1,\dots ,\alpha (x_{g_n})=g_n$ of the factors $x_{g_1},
\dots , x_{g_n}$ and denote by $G_O$ the subgroup of $G$ generated
by the elements of $O$.

\begin{claim} {\rm (\cite{Ku1})} \label{cc} The subgroup $G_s$ of $G$ is well defined, that is,
it does not depend on a presentation of $s$ as a product of
generators $x_{g_i}\in \mathcal X_O$.
\end{claim}

For subgroups $H$ and $K$ of a group $G$, we put
$$S(G,O)^H=\{s\in S(G,O)\mid G_s=H\},$$
$$S(G,O)_{K}= \{ s\in S(G,O)\mid \alpha (s)\in K \},$$ and
$S(G,O)_{K}^H= S(G,O)_{K}\cap S(G,O)^H$. It is easy to see that
$S(G,O)^H$ (respectively $S(G,O)_{K}^H$) is isomorphic to the
semigroup $S(H,H\cap O)^H$ (respectively, isomorphic to $S(H,H\cap
O)_{K\cap H}^H$) and the isomorphism is induced by the embedding
$(H,H\cap O)\hookrightarrow (G,O)$.

\begin{prop} {\rm (\cite{Ku1})} \label{simple} Let $(G,O)$ be an equipped group and let $s\in S(G,O)$. We have
\begin{itemize} \item[($1$)] $\text{ker}\, \rho $ coincides with the centralizer
$C_{O}$ of the group $G_{O}$ in $G$;
\item[($2$)] if
$\alpha (s)$ belongs to the center $Z(G_s)$ of $G_s$, then for each
$g\in G_s$ the action $\rho (g)$ leaves fixed the element $s\in
S(G,O)$;
\item[($3$)] if $\alpha (s\cdot x_g)$ belongs to the center $Z(G_{s\cdot x_g})$ of
$G_{s\cdot x_g}$, then $s\cdot x_g=x_g\cdot s$,
\item[($4$)] if $\alpha (s)=\bf{1}$, then $s\cdot s'=s'\cdot s$
for any $s'\in S(G,O)$.
\end{itemize}
\end{prop}

\begin{claim}{\rm (\cite{Ku1})} \label{commut} For any equipped group $(G,O)$ the
semigroup $S(G,O)_{\bf{1}}$ is contained in the center of the
semigroup $S(G,O)$ and, in particular, it is commutative.
\end{claim}

Note that if $g\in O$ is an element of order $n$, then $x_g^n\in
S(G,O)_{\bf{1}}$.
\begin{lem} {\rm (\cite{Ku1})} \label{fried}  Let $s\in
S(G,O)_{Z(G_O)}$ and $s_1\in S(G,O)^{G_O}$, where $Z(G_O)$ is the
center of $G_O$. Then
\begin{equation} \label{oo} s\cdot s_1= \rho(g)(s)\cdot s_1
\end{equation}
for all $g\in G_O$.

In particular, if $s\in S(G,O)^G$, $C\subset O$ is a conjugacy class
of $G$, and $g_1^n$ belongs to the center $Z(G)$ of $G$ for certain
$g_1\in C$, then for any $g_2\in C$ we have
\begin{equation} \label{pp} x_{g_1}^{n}\cdot s=x_{g_2}^{n}\cdot s.
\end{equation}
\end{lem}
\begin{prop} {\rm (\cite{Ku1})} \label{simpl}
The elements of $S(G,O)_{\bf{1}}^G$ are fixed under the conjugation
action of $G$.
\end{prop}

\subsection{Factorization semigroups over equivalent equipped
groups} \label{sub3} Here, we introduce an additional assumption
with regard to $O$ in an equipped group $(G,O)$: we assume that
\begin{enumerate}
\item[(iv)]
the elements of $O$ generate the group $G$.
\end{enumerate}

In \cite{Ku3}, a $C$-graph was associated with each equipped group.
By definition, the $C$-{\it graph $\Gamma_{(G,O)}$ of an equipped
group $(G,O)$} is a directed labeled graph. Its vertices are labeled
by elements of $O$, and this labeling is a bijection between $O$ and
the set of vertices, $V=\{ v_{g}\mid g\in O\}$. Each edge of
$\Gamma_{(G,O)}$ also is labeled by an element of $O$. Namely, two
vertices $v_{g_1}$ and $v_{g_2}$, $g_1, g_2\in O$, are connected by
an edge $e_{v_{g_1},v_{g_2},g}$ with label $g\in O$ if and only if
$g^{-1}g_1g=g_2$. (A $C$-graph may contain loops, and several edges,
but with distinct labels, may  connect the same pair of vertices in
the same direction.)

Obviously, the conjugacy classes $C_i\subset O$, $1\le i\le m$, are
in one-to-one correspondence with the connected components
$\Gamma_i$ of the $C$-graph $\Gamma_{(G,O)}$; more precisely,
$v_{g}\in \Gamma_i$ if and only if $g\in C_i$.

Two equipped groups $(G_1,O_1)$ and $(G_2,O_2)$ are called {\it
equivalent} if their $C$-graphs, $\Gamma_{(G_1,O_1)}$ and
$\Gamma_{(G_2,O_2)}$, are isomorphic as $C$-graphs; in other words,
if there is a bijection $O_1 \to O_2$ that induces an isomorphism of
$C$-graphs between $\Gamma_{(G_1,O_1)}$ and $\Gamma_{(G_2,O_2)}$.

To each $C$-graph $\Gamma=\Gamma_{(G,O)}$ one associates  a
$C$-group $G_{\Gamma}=(\widetilde G,\widetilde O)$ equivalent to
$(G,O)$. Denoting by $g\mapsto \tilde g$ the bijection $O\to \tilde
O$, we can describe $\tilde G$ as the group defined by generators
$\tilde g\in\tilde O$ and the relations
$$
\{\widetilde g_{3}^{-1}\widetilde g_{1}\widetilde g_{3}=\widetilde
g_{2}\,\text{ if  and only if there is an edge} \,
e_{v_{g_1},v_{g_2},g_3}\in\Gamma\}
$$
(equivalently, $\widetilde g_{3}^{-1}\widetilde g_{1}\widetilde
g_{3}=\widetilde g_{2}$ if and only if the relation
$g_{3}^{-1}g_{1}g_{3}=g_{2}$ holds in $G$). These generators are
called {\it $C$-generators}. They are in one-to-one correspondence
with the vertices of $\Gamma$ due to the composed bijection $\tilde
g\mapsto g\mapsto v_g$. Furhermore, two $C$-generators $\widetilde
g_1$ and $\widetilde g_2$ belong to the same connected component of
$\Gamma_{(G,O)}$ if and only if they are conjugate.
The set
$\widetilde O$ of $C$-generators of the $C$-group $\widetilde G$
satisfies the assumptions
(i)-(iv).

The one-to-one map $\beta_{(G,O)}: \widetilde O\to O$ given by
$\beta_{(G,O)}(\widetilde g)=g$ defines an epimorphism
$\beta=\beta_{(G,O)}: (\widetilde G,\widetilde O)\to (G,O)$ of
equipped groups and an isomorphism $\beta_*: S(\widetilde
G,\widetilde O)\to  S(G,O)$ of semigroups. By Claim 8 in \cite{Ku3},
$\ker \beta$ is a subgroup of the center $Z(\widetilde G)$ of the
$C$-group $\widetilde G$.

By adding the commutativity relations one shows that the
abelianization $H_1(\widetilde G,\mathbb Z)=\widetilde G/[\widetilde
G,\widetilde G]$ of $\widetilde G$ is isomorphic to $\mathbb Z^m$.
Moreover, due to a fixed ordering of the conjugacy classes $\{
\widetilde C_1,\dots, \widetilde C_m\}$ it comes with a natural
basis, where in terms of the abelianization homomorphism $\text{ab}:
\widetilde G\to H_1(\widetilde G,\mathbb Z)$ the $i$-th element of
the basis is given by $\text{ab}(\widetilde g)$ with $g\in C_i$,
$1\le i\le m$.

The homomorphism $\tau=\text{ab}\, \circ\, \beta_*^{-1}:S(G,O)\to \mathbb
Z^m_{\geq 0}$ is called the {\it type} homomorphism, the image
$\tau(s)=(\tau_1(s),\dots,\tau_m(s))\in \mathbb Z^m_{\geq 0}$ is
called the {\it type} of $s\in S(G,O)$, and the $i$th coordinate
$\tau_i(s)$ of $\tau(s)$ is called the {\it $i$th type} of $s$.

A $C$-group $\widetilde G$ is called {\it $C$-finite} if the number
of vertices of the graph $\Gamma_{(\widetilde G,\widetilde O)}$ is
finite. By Proposition 3 in \cite{Ku3}, the commutant $[\widetilde
G,\widetilde G]$ of a $C$-finite group $\widetilde G$ is a finite
group. The order $a_{(G,O)}=|\ker \beta \cap [\widetilde
G,\widetilde G]|$ of the group $\ker \beta \cap [\widetilde
G,\widetilde G]$ is called the {\it ambiguity index} of an equipped
finite group $(G,O)$. If $O^{\prime}\subset O$ are two  equipments
of a finite group $G$ such that the elements of $O^{\prime}$
generate the group $G$, then by Corollary 2 in \cite{Ku3}, we have
$a_{(G,O)}\leq a_{(G,O')}$.

In \cite {Ku3} the following theorems are proved.

\begin{thm} \label{Cuniq} Let $(G,O)$, $O=C_1\sqcup \dots \sqcup
C_m$, be an equipped finite group and $(\widetilde G,\widetilde
O)=G_{\Gamma}$ with $\Gamma=\Gamma_{(G,O)}$ the $C$-group equivalent
to $(G,O)$. Then there is a constant $T\in \mathbb N$ such that for
any element $s_1\in S(G,O)^G$ with $\tau_i(s_1)\geq T$ for all
$i=1,\dots, m$ there exist $a_{(G,O)}$ elements $s_1,\dots,
s_{a_{(G,O)}}\in S(G,O)^G$ such that
\begin{itemize}
\item[$(1)$] $s_i\neq s_j$ for $1\leq i< j\leq a_{(G,O)}$;
\item[$(2)$] $\tau(s_i)=\tau(s_1)$ for $1\leq i\leq a_{(G,O)}$;
\item[$(3)$] $\alpha_G(s_i)=\alpha_G(s_1)$ for $1\leq i\leq a_{(G,O)}$;
\item[$(4)$] if $s\in S(G,O)^G$,
$\tau(s)=\tau(s_1)$ and $\alpha_G(s)=\alpha_G(s_1)$, then $s=s_i$
for some $i$, $1\leq i\leq  a_{(G,O)}$.
\item[$(5)$] if $s\in S(G,O)^G$ and $\alpha_{\widetilde G}(s)=\alpha_{\widetilde G}(s_1)$, then $s=s_1$.
\end{itemize}
\end{thm}

\begin{thm} \label{Cuniq2} Let $G$ be a finite group and
$O^{\prime}\subset O$ be two its equipments such that the elements
of $O^{\prime}=C_1\sqcup \dots \sqcup C_k$ generate the group $G$.
Then there is a constant $T=T_{O^{\prime}}$ such that if for an
element $s_1\in S(G,O)^G$ the $i$th type $\tau_i(s_1)\geq T$ for all
$i=1,\dots, k$, then there are not more than $a_{(G,O^{\prime})}$
elements $s_1,\dots, s_{n}\in S(G,O)^G$ such that
\begin{itemize}
\item[$(i)$] $s_i\neq s_j$ for $1\leq i< j\leq n$;
\item[$(ii)$] $\tau(s_i)=\tau(s_1)$ for $1\leq i\leq n$;
\item[$(iii)$] $\alpha_G(s_i)=\alpha_G(s_1)$ for $1\leq i\leq n$,
\end{itemize}
where $a_{(G,O^{\prime})}$ is the ambiguity index of
$(G,O^{\prime})$.
\end{thm}

Theorem \ref{Cuniq2} is exactly Theorem 7 from \cite{Ku3}. Theorem
\ref{Cuniq}, items (1)--(4), is Theorem 6 from \cite{Ku3}, while the
item (5) is a direct consequence of (\cite{Ku3}, Theorems 5 and 6)
and the following straightforward remark.

\begin{rem} \label{rem} {\rm The elements $s_1,\dots, s_{a_{(G,O)}}\in S(G,O)^G$,
whose existence is claimed by Theorem \ref{Cuniq}, are distinguished
by their valuers $\alpha_{\widetilde G}(s_i) \in \widetilde G$.
Namely, for $i\neq j$ the element $\alpha_{\widetilde
G}(s_i)\alpha_{\widetilde G}(s_j)^{-1}$ is a non-trivial element of
$\ker \beta_{(G,O)}\cap [\widetilde G,\widetilde G]$.}
\end{rem}

\subsection{Admissible subgroups of
$Aut(\mathbb F^{n+2p})$} \label{adm} In this subsection, in order to
introduce a notion of algebraic covering semigroups, we pick out
some class of subgroups of the automorphism groups of free groups,
called {\it admissible automorphism groups.}

Let $\mathbb F^{n+2p}$ be a free group freely generated by $n+2p$
elements. Let $\mathcal G=\{ \gamma_1,\dots, \gamma_n\}\subset
\mathbb F^{n+2p}$, $\mathcal L=\{ \lambda_1,\dots,\lambda_p\}\subset
\mathbb F^{n+2p}$, and $\mathcal M=\{ \mu_1,\dots,\mu_p\}\subset
\mathbb F^{n+2p}$ be three ordered sets such that the elements of
$\mathcal B=\mathcal G\cup \mathcal L\cup\mathcal M$ generate the
group $\mathbb F^{n+2p}$.

Let $\overline n=(n_1,\dots, n_{p+1})$ be an ordered non-negative
partition of the number $n$, that is, an ordered $(p+1)$-tuple of
non-negative integers whose sum is equal to $n$:
$$n=n_1+\dots +n_{p+1}, \qquad n_i\in \mathbb Z_{\geq 0}.$$
We put $k_i=\sum_{j=1}^in_j$. Each partition $\overline n$ defines
its own ordering on $\mathcal B$:
$$\gamma_1,\dots, \gamma_{k_1},\lambda_1,\mu_1, \dots,
\gamma_{k_{i-1}+1},\dots,\gamma_{k_{i}},\lambda_i,\mu_i,\gamma_{k_{i}+1},\dots,\gamma_{k_{i+1}},
\dots , \lambda_p,\mu_p,\gamma_{k_{p}+1}\dots,\gamma_{n}$$ (here the
set $\{\gamma_{k_{i-1}+1},\dots,\gamma_{k_{i}}\}$ is empty if
$n_i=0$). Denote by $\mathcal B_{\overline n}$ the set $\mathcal B$
with the ordering defined by partition $\overline n$ and call it a
{\it frame} of  $\mathbb F^{n+2p}$. In particular, $\mathcal
B_{(n,0,\dots,0)}=\{
\gamma_1,\dots,\gamma_n,\lambda_1,\mu_1,\dots,\lambda_p,\mu_p\}$.
The element
$$\partial \mathcal
B_{\overline n}=\gamma_1\dots \gamma_{k_1}[\lambda_1,\mu_1] \dots
\gamma_{k_{i-1}+1}\dots\gamma_{k_{i}}[\lambda_i,\mu_i]\gamma_{k_{i}+1}\dots\gamma_{k_{i+1}}
\dots [\lambda_p,\mu_p]\gamma_{k_{p}+1}\dots\gamma_{n}$$ of $\mathbb
F^{n+2p}$ is called {\it the boundary of} $\mathcal B_{\overline
n}$.

Given a set $\mathcal B'=\mathcal G'\cup \mathcal L'\cup\mathcal M'$
as above and two adjacent partitions, $\overline n'=(\dots,
n_{i-1},n_i,n_{i+1},n_{i+2},\dots)$  and $\overline n''=(\dots,
n_{i-1},n_i-1,n_{i+1}+1,n_{i+2},\dots)$, we define {\it an
elementary frame change} $h_{\overline n',\overline n''}$ that
results both in change of the generating set and the ordering.
Namely, we put
$$
h_{\overline n',\overline n''}(\mathcal B'_{\overline n'})=\mathcal
B''_{\overline n''},
$$
where $\mathcal B''=\mathcal G''\cup \mathcal L''\cup\mathcal M''$
with $\lambda''_j=\lambda'_j$ and $\mu''_j=\mu'_j$ for all
$j=1,\dots, p$, while $\gamma''_j=\gamma'_j$ for $j\neq
k_i=n_1+\dots+n_i$ and
$\gamma''_{k_i}=([\lambda'_{i},\mu'_{i}])^{-1}\gamma'_{k_i}[\lambda'_{i},\mu'_{i}]$.
The inverse change $h^{-1}_{\overline n',\overline n''}=h_{\overline
n'',\overline n'}$ also will be called an elementary frame change.
Two frames of $\mathbb F^{n+2p}$ are said to be {\it strongly
equivalent} if one of them can be obtained from the other one by a
finite sequence of elementary frame changes. Any composition of
elementary changes transforming a frame $\mathcal B'_{\overline n'}$
into a frame $\mathcal B''_{\overline n''}$ will also be denoted by
$h_{\overline n',\overline n''}$.

The proof of the following properties is straightforward.
\begin{claim} \label{change} Let $\mathcal B'_{\overline n}$ and $\mathcal B''_{\overline n}$
be two frames strongly equivalent to a frame $\mathcal
B_{(n,0,\dots,0)}$. Then $\mathcal B'_{\overline n}=\mathcal
B''_{\overline n}$.
\end{claim}
\begin{claim} \label{change} Let $\mathcal B'_{\overline n'}$ and $\mathcal B''_{\overline n''}$
be two strongly equivalent frames. Then  $\partial \mathcal
B'_{\overline n'}=\partial \mathcal B''_{\overline n''}$. \qed
\end{claim}

The group $Aut(\mathbb F^{n+2p})$ naturally acts on the set of
frames. This action respects the partition. Given $h\in Aut(\mathbb
F^{n+2p})$ and a frame $\mathcal B_{\overline n}$, we put $\mathcal
B'=h(\mathcal B)$ and $h(\mathcal B_{\overline n})=\mathcal
B'_{\overline n}$. As usually, the orbit of $\mathcal B_{\overline
n}$ under the action of a subgroup $H$ of $Aut(\mathbb F^{n+2p})$ is
denoted by $H\mathcal B_{\overline n}$. The following Lemma is
obvious.

\begin{lem} \label{import} Let $H$ be a subgroup of $Aut(\mathbb F^{n+2p})$ and $\mathcal
B'_{\overline n'}$, $\mathcal B''_{\overline n''}$ two strongly
equivalent frames. Then
\begin{itemize}
\item[($i$)]
for any $h\in H$, the frames $h(\mathcal B'_{\overline n'})$ and
$h(\mathcal B''_{\overline n''})$ are strongly equivalent and
$h(\mathcal B'_{\overline n'})= h_{\overline n',\overline
n''}^{-1}(h(B''_{\overline n''}))$,
\item[($ii$)] the map
$h_{\overline n',\overline n''}:H\mathcal B'_{\overline n'}\to
H\mathcal B''_{\overline n''}$ is one-to-one.
\end{itemize}
\end{lem}

Let us fix a frame $\mathcal B_{{\bf{1}}}=\mathcal
B_{(n,0,\dots,0)}=\{
\gamma_1,\dots,\gamma_n,\lambda_1,\mu_1,\dots,\lambda_p,\mu_p\}$
and, for each $i$ with $2\le i\le p+1$, put
$\mathcal B_{{\bf{i}}}=h_{\overline n, \overline n'}\mathcal B_{{\bf{1}}}$
where $\overline n=(n,0,\dots,0)$ and $\overline
n'=(n-1,0,\dots,0,1,0,\dots,0)$ with $1$ on the $i$-th place.

We specify several auxiliary automorphisms of $\mathbb F^{n+2p}$.

The automorphism $\sigma_i$ with $i=1,\dots, n-1$ is defined by its
action in the frame $\mathcal B_{{\bf{1}}}$ as follows:
$$\begin{array}{lll} \sigma_i(\lambda_j) & =\lambda_j \qquad \text{for}\, j=1,\dots, p, \\
\sigma_i(\mu_j) & =\mu_j \qquad \text{for}\, j=1,\dots, p, \\
\sigma_i(\gamma_j) & =\gamma_j \qquad  \text{for}\, j\neq i,i+1,
\\ \sigma_i(\gamma_i) & =\gamma_{i+1}, & \\
\sigma_i(\gamma_{i+1}) & =\gamma_i^{\gamma_{i+1}}. &
\end{array} $$

The automorphism $\xi_{i,\lambda}$ with $i=1,\dots, p$ is defined by
its action in the frame $\mathcal B_{{\bf{i}}}$ as follows:
$$\begin{array}{lll} \xi_{i,\lambda}(\lambda_{j,{\bf{i}}}) & =\lambda_{j,{\bf{i}}} \qquad \text{for}\, j\neq i, \\
\xi_{i,\lambda}(\mu_{j,{\bf{i}}}) & =\mu_{j,{\bf{i}}} \qquad \text{for}\, j=1,\dots, p, \\
\xi_{i,\lambda}(\gamma_{j,{\bf{i}}}) & =\gamma_{j,{\bf{i}}} \qquad
\text{for}\, j\neq n,
\\ \xi_{i,\lambda}(\gamma_{n,{\bf{i}}}) & =\gamma_{n,{\bf{i}}}^{c_{1,i}}, \,\,\,\,
\,\, \text{where}\,\,
c_{1,i}=\lambda_{{i,{\bf{i}}}}\mu_{{i,{\bf{i}}}}^{-1}\lambda_{{i,{\bf{i}}}}^{-1}\gamma_{n,{\bf{i}}}^{-1}, \\
\xi_{i,\lambda}(\lambda_{{i,{\bf{i}}}}) &
=\gamma_{n,{\bf{i}}}\lambda_{{i,{\bf{i}}}}. &
\end{array} $$

The automorphism $\xi_{i,\mu}$ with $i=1,\dots, p$ is defined by its
action in the frame $\mathcal B_{{\bf{i}}}$ as follows:
$$\begin{array}{lll} \xi_{i,\mu}(\mu_{j,{\bf{i}}}) & =\mu_{j,{\bf{i}}} \qquad \text{for}\, j\neq i, \\
\xi_{i,\mu}(\lambda_{j,{\bf{i}}}) & =\lambda_{j,{\bf{i}}} \qquad \text{for}\, j=1,\dots, p, \\
\xi_{i,\mu}(\gamma_{j,{\bf{i}}}) & =\gamma_{j,{\bf{i}}} \qquad
\text{for}\, j\neq n,
\\ \xi_{i,\mu}(\gamma_{n,{\bf{i}}}) & =\gamma_{n,{\bf{i}}}^{c_{2,i}},
\,\,\,\, \, \,   \text{where}\,\,
c_{2,i}=\mu_{{i,{\bf{i}}}}\lambda^{-1}_{{i,{\bf{i}}}}\mu^{-1}_{{i,{\bf{i}}}}\gamma_{n,{\bf{i}}}, \\
\xi_{i,\mu}(\mu_{{i,{\bf{i}}}}) &
=\gamma_{n,{\bf{i}}}^{-1}\mu_{{i,{\bf{i}}}}. &
\end{array} $$

The automorphism $\zeta_{i}$ with $i=1,\dots, p$ is defined by its
action in the frame $\mathcal B_{{\bf{i}}}$ as follows:
$$\begin{array}{ll} \zeta_i(\lambda_{j,{\bf{i}}}) & =\lambda_{j,{\bf{i}}} \qquad \text{for}\, j\neq i,
\\ \zeta_i(\mu_{j,{\bf{i}}}) & =\mu_{j,{\bf{i}}} \qquad \text{for}\, j\neq i,
\\ \zeta_i(\gamma_{j,{\bf{i}}}) & =\gamma_{j,{\bf{i}}} \qquad
\text{for}\, j\neq n,
\\ \zeta_i(\gamma_{n,{\bf{i}}}) & =c_{3,{\bf{i}}}, 
\\ \zeta_i(\lambda_{i,{\bf{i}}})
& =\lambda_{i,{\bf{i}}}^{c_{3,{\bf{i}}}},
\\ \zeta_i(\mu_{i,{\bf{i}}}) &
=\mu_{i,{\bf{i}}}^{c_{3,{\bf{i}}}},
\end{array} $$
where
$c_{3,{\bf{i}}}=\gamma_{n,{\bf{i}}}^{[\lambda_{i,{\bf{i}}},\mu_{i,{\bf{i}}}]}$.

For $0\leq p_1\leq p$, denote by $Br_{n,p_1}$ the subgroup of the
group $Aut(\mathbb F^{n+2p})$ generated by the elements
$\sigma_1,\dots,\sigma_{n-1},\xi_{1,\lambda},\dots
,\xi_{p_1,\lambda},\xi_{1,\mu},\dots,
\xi_{p_1,\mu},\zeta_{1},\dots,\zeta_{p_1}$. Obviously, for $p_1\leq
p$ the group $Br_{n,p_1}$ is a subgroup of $Br_{n,p}$ and the groups
$Br_{n,p_1}\subset Aut(\mathbb F^{n+2p_1})$ and $Br_{n,p_1}\subset
Br_{n,p}\subset Aut(\mathbb F^{n+2p})$ are naturally isomorphic. The
groups $Br_{n,p}$ will be called {\it algebraic braid groups.}

\begin{claim} \label{fixed} The boundary
$\partial \mathcal
B_{{\bf{1}}}=\gamma_1\dots\gamma_n[\lambda_1,\mu_1]\dots[\lambda_p,\mu_p]\in
\mathbb F^{n+2p}$ is fixed under the action of $Br_{n,p}$.
\end{claim}
\proof Obviously, $\partial \mathcal B_{{\bf{1}}}$ is fixed under
the actions of $\sigma_i$ for $i=1,\dots,n-1$ and the actions of
$\xi_{1,\lambda}$, $\xi_{1,\mu}$, $\zeta_{1}$, as well as $\partial
\mathcal B_{{\bf{i}}}$ with $i\geq 2$ is fixed under the actions of
the automorphisms $\xi_{i,\lambda}$, $\xi_{i,\mu}$, $\zeta_{i}$.
Now, the statement follows from Claim \ref{change}.
\qed \\

Let $(\dots,\gamma'_i,\gamma'_{i+1},\dots)$ be a part of a frame
$\mathcal B'_{\overline n}$ that we assume to be strongly equivalent
to $\mathcal B_{{\bf{1}}}$. Denote by $\sigma_{i,\overline n}$ an
automorphism of $\mathbb F^{n+2p}$ such that $\sigma_{i,\overline
n}(\gamma'_i)=\gamma'_{i+1}$, $\sigma_{i,\overline
n}(\gamma'_{i+1})=(\gamma'_i)^{\gamma'_{i+1}}$, and
$\sigma_{i,\overline n}$ leaves fixed all the other elements of
$\mathcal B'_{\overline n}$. Similarly, if $(\dots,\gamma'_j,
\lambda'_i,\mu'_i,\dots ...)$ is a part of a frame $\mathcal
B'_{\overline n}$, then denote by $\xi_{i,\overline n,\lambda}$,
$\xi_{i,\overline n,\mu}$, and  $\zeta_{\overline n,i}$ the
automorphisms of $\mathbb F^{n+2p}$ that leave fixed all the
elements of $\mathcal B'_{\overline n}$ except $\gamma'_j$,
$\lambda'_i$, $\mu'_i$ and act on $\gamma'_j$, $\lambda'_i$,
$\mu'_i$ by the same formulas as $\xi_{i,\lambda}$, $\xi_{i,\mu}$,
and $\zeta_{\overline n,i}$ act on the elements $\gamma_n$,
$\lambda_i$, $\mu_i$ of the frame $\mathcal B_{{\bf{i}}}$ (we just
replace $n$ by $j$).

\begin{lem}\label{imp} Let $\mathcal B'_{\overline n}$ be strongly equivalent to $\mathcal
B_{{\bf{1}}}$. Then the automorphisms $\sigma_{i,\overline n}$,
$\xi_{i,\overline n,\lambda}$, $\xi_{i,\overline n,\mu}$, and
$\zeta_{\overline n,i}$ of $\mathbb F^{n+2p}$ belong to $Br_{n,p}$.
\end{lem}

\proof Follows from Claim \ref{change} and Lemma
\ref{import} by straightforward induction on $p$. \qed \\

We say that a subgroup $H_{n,p}$ of $Aut(\mathbb F^{n+2p})$ is {\it
admissible} if:
\begin{itemize}
\item[$1$)] $Br_{n,p}\subseteq H_{n,p}$;
\item[$2$)] for each $h\in
H_{n,p}$ there is a permutation  $\sigma_h \in \mathcal S_n$ such
that $h(\gamma_i)$ is conjugate to $\gamma_{\sigma_h(i)}$;
\item[$3$)] for each $h\in
H_{n,p}$ it holds the relation
$$
h(\gamma_1\dots
\gamma_{n}[\lambda_1,\mu_1]\dots[\lambda_{p},\mu_{p}])=\gamma_1\dots
\gamma_{n}[\lambda_1,\mu_1]\dots[\lambda_{p},\mu_{p}]$$ (here
$\gamma_1\dots \gamma_{n}={\bf{1}}$ if $n=0$).
\end{itemize}

Let us fix a frame $\mathcal B_{{\bf{1}}}=\{
\gamma_1,\dots,\gamma_n,\lambda_1,\mu_1,\dots,\lambda_p,\mu_p\}$ of
$\mathbb F^{n+2p}$ and let $f: \mathbb F^{n+2p}\to G$ be a
homomorphism to an equipped group $(G,O)$ such that $f(\gamma_i)\in
O$ (we call such an $f$ {\it an equipped homomorphism} to $(G,O)$).
Put $g_i=f(\gamma_i)$ for $1\le i\le n$ and $a_j=f(\lambda_j)$,
$b_j=f(\mu_j)$ for $1\le j\le p$.

To each frame $\mathcal B_{\overline n}$ strongly equivalent to
$\mathcal B_{{\bf{1}}}$, we associate a word $W_{f,\mathcal
B_{\overline n}}$ in the alphabet $\mathcal Z=\mathcal Z_{(G,O)}=
\mathcal X_O\cup \mathcal Y_G$, where $\mathcal X_O=\{ x_g\mid g\in
O\}$ is the alphabet we used already in subsection \ref{semigr} and
$\mathcal Y_G=\{ y_{a,b}\mid (a,b)\in G^2\}$. We put
$$
W_{f,\mathcal B_{{\bf{1}}}}=x_{g_1}\dots x_{g_n}y_{a_1,b_1}\dots
y_{a_p,b_p}
$$
and then construct the words $W_{f,\mathcal B_{\overline n}}$
iteratively by {\it elementary moves}: in notation used in the
definition of an elementary frame change $ h_{\overline n',\overline
n''}(\mathcal B'_{\overline n'})=\mathcal B''_{\overline n''}$,
where $\overline n'=(\dots, n_{i-1},n_i,n_{i+1},n_{i+2},\dots)$  and
$\overline n''=(\dots, n_{i-1},n_i-1,n_{i+1}+1,n_{i+2},\dots)$ are
two adjacent partitions, the elementary move $W_{f,\mathcal
B'_{\overline n'}}\mapsto W_{f,\mathcal B''_{\overline n''}}$
consists in the replacement of two adjacent letters
$x_{g'_{k_i}}y_{a'_i,b'_i}$  in $W_{f,\mathcal B'_{\overline n'}}$
by $y_{a'_i,b'_i}x_{([a_i,b_i])^{-1}g'_{k_i}[a'_i,b'_i]}$ (as in the
definition of elementary frame changes, $k_i=n_1+\dots+n_i$).

Denote by $\overline W_f(G,O)$ the set of words which can be
obtained from $W_{f,\mathcal B_{{\bf{1}}}}$ by finite sequences of
elementary moves and put $\overline W_{n,p}(G,O)=\bigcup_f\overline
W_f(G,O)$, where the union is taken over all equipped homomorphisms
$f:\mathbb F^{n+2p}\to (G,O)$. We say that two words are {\it
$qf$-equivalent} if they belong to the same set $\overline
W_f(G,O)$.

Every admissible group $H_{n,p}$ acts on $\overline W_{n,p}(G,O)$.
Namely, we put
$$
h(W_{f,\mathcal B_{\overline n}})=W_{f,\mathcal B'_{\overline n}},
B'_{\overline n}= h(B_{\overline n}).
$$
In particular, if $W_{f,\mathcal B_{\overline n}}$ is obtained from
$W_{f,\mathcal B_{{\bf{1}}}}$ by a finite sequence of elementary
moves, then $h(W_{f,\mathcal B_{\overline n}})$ is obtained from
$h(W_{f,\mathcal B_{{\bf{1}}}})$ by the same sequence of elementary
moves.

Let $h$ be an element of an admissible group $H_{n,p}$. We have
$h(\gamma_i)=\gamma_{\sigma_h(i)}^{w_i}$, $h(\lambda_i)=u_i$, and
$h(\mu_i)=v_i$, where $w_i$, $u_i$, $v_i$ are some elements of
$\mathbb F^{n+2p}$. Denote by the same letters the words
$w_i=w_i(\gamma_1,\dots, \gamma_n, \lambda_1,\mu_1,\dots,
\lambda_p,\mu_p)$, $u_i=u_i(\gamma_1,\dots, \gamma_n,
\lambda_1,\mu_1,\dots, \lambda_p,\mu_p)$, and
$v_i=v_i(\gamma_1,\dots, \gamma_n, \lambda_1,\mu_1,\dots,
\lambda_p,\mu_p)$ in letters $\gamma_1,\dots, \gamma_n,
\lambda_1,\mu_1,\dots, \lambda_p,\mu_p$ and their inverses
representing these elements in $\mathbb F^{n+2p}$. Consider elements
$g_1,\dots ,g_n$, $a_1,b_1,\dots,a_p,b_p$ of an equipped group
$(G,O)$, where $g_1,\dots,g_n\in O$, and let us substitute $g_j$ for
$\gamma_j$, $a_j$ for $\lambda_j$, $b_j$ for $\mu_j$ into the words
$w_i$, $u_i$, $v_i$ and denote the corresponding elements of $G$ by
$\overline w_i=w_i(g_1,\dots g_n,a_1,b_1,\dots,a_p,b_p)$, $\overline
u_i=u_i(g_1,\dots g_n,a_1,b_1,\dots,a_p,b_p)$, and $\overline
v_i=v_i(g_1,\dots g_n,a_1,b_1,\dots,a_p,b_p)$. Denote by $\langle
g_1, \dots, g_n,a_1,b_1,\dots,a_p,b_p\rangle$ a subgroup of $G$
generated by the elements $g_1,\dots, g_n,a_1,b_1,\dots,a_p,b_p\in
G$.

\begin{claim} \label{ggg} In notations and assumptions used above, we have
$$\langle g_1,\dots g_n,a_1,b_1,\dots ,a_p,b_p\rangle= \langle
g_{\sigma_h(1)}^{\overline w_1},\dots, g_{\sigma_h(n)}^{\overline
w_n}, \overline u_1,\overline v_1,\cdots, \overline u_p,\overline
v_p \rangle .$$
\end{claim}
\proof It suffices to note that the subgroup $\langle g_1,\dots
g_n,a_1,b_1,\dots ,a_p,b_p\rangle$ of $G$ is the image in $G$ of the
group $\mathbb F^{n+2p}$ under the homomorphism $f:\mathbb
F^{n+2p}\to G$ given by $f(\gamma_i)=g_i$, $f(\lambda_j)=a_j$, and
$f(\mu_j)=b_j$. \qed

\subsection{Definition of covering semigroups} \label{def-cov} Let
$(G,O)$ be an equipped group. Denote by  $F\mathbb S(G,O)$ the free
semigroup over the alphabet $\mathcal Z=\mathcal Z_{(G,O)}= \mathcal
X_O\cup \mathcal Y_G$ introduced in subsection \ref{adm} and call
$F\mathbb S(G,O)$ {\it free covering semigroup over the equipped
group} $(G,O)$.

All the covering semigroups considered below are factor semigroups
of $F\mathbb S(G,O)$. In particular, this is the case of what we
call the {\it quasi-free algebraic covering semigroup} $qF\mathbb
S(G,O)$ that we define as a semigroup generated by the alphabet
$\mathcal Z$ and subject to relations
\begin{equation} \label{rel12'} x_{g}\cdot
y_{a,b}=y_{a,b}\cdot x_{g^{[a,b]}},\quad g\in O,\, \, \,  a,b\in G 
\end{equation}
(in other words, the elements of $qF\mathbb S(G,O)$ are the sets of
$qf$-equivalent words (see subsection \ref{adm})).

We follow notation of Subsection \ref{adm}. Let $\mathcal H=\{
H_{n,p}\}_{\{ n\geq 0,p\geq 0\}}$ be a collection of automorphism
groups that satisfy conditions $2), 3)$ from the definition of
admissible automorphism groups. We associate with each $h\in
H_{n,p}$ a set $R_h$ of relations
\begin{equation} \label{admis}
x_{g_1}\cdot .\, .\, . \cdot x_{g_n}\cdot y_{a_1,b_1}\cdot .\, .\,
.\cdot y_{a_p,b_p}= x_{g_{\sigma_h(1)}^{\overline w_1}}\cdot .\, .\,
.\cdot x_{g_{\sigma_h(n)}^{\overline w_n}}\cdot y_{\overline
u_1,\overline v_1}\cdot .\, .\, .\, . \cdot y_{\overline
u_p,\overline v_p}\end{equation} taken over all $(g_1,\dots, g_n)\in
O^n$ and all $(a_1,b_1,\dots,a_p,b_p)\in G^{2p}$. Denote by
$$\displaystyle \mathcal R_{\mathcal H}= \bigcup_{n,p}(\bigcup_{h\in
H_{n,p}}R_h)
$$
and  consider a factor semigroup $qF\mathbb S(G,O)/ {\mathcal
R_{\mathcal H}}$. In particular, the semigroup
$$\mathbb S(G,O)=qF\mathbb S(G,O)/\{ \mathcal R_{\mathcal B}\}$$ is called
{\it the strong covering semigroup}, where $\mathcal B=\{
Br_{n,p}\}_{\{ n\geq 0,p\geq 0\}}$. For a collection  $\mathcal H=\{
H_{n,p}\}_{\{ n\geq 0,p\geq 0\}}$  of admissible automorphism
groups,  we denote the semigroup $qF\mathbb S(G,O)/ {\mathcal
R_{\mathcal H}}$  by $\mathbb S_{\mathcal H-\text{equiv}}(G,O)$ and
call it an {\it admissible covering semigroup}.

\begin{prop} \label{red-rel} The strong covering semigroup $\mathbb S(G,O)$
is isomorphic to the semigroup generated by the alphabet $\mathcal
Z=\mathcal Z_{(G,O)}=\mathcal X_O\cup \mathcal Y_G$ and subject to
relations
\begin{equation} \label{rel11} x_{g_1}\cdot
x_{g_2}=x_{g_2}\cdot\,x_{g_1^{g_2}}\end{equation} for any $x_{g_1},
x_{g_2}\in \mathcal X_O$, and
\begin{equation} \label{rel12} x_{g}\cdot
y_{a,b}=y_{a,b}\cdot x_{g^{[a,b]}},\end{equation} \begin{equation}
\label{rel13} x_{g}\cdot y_{a,b}=x_{g^{c_1}}\cdot y_{ga,b}, \qquad
c_1=ab^{-1}a^{-1}g^{-1},\end{equation}
\begin{equation} \label{rel16} y_{a,b} \cdot x_{g}=y_{a,g^{-1}b}\cdot x_{g^{c_2}}, \qquad c_2=ba^{-1}b^{-1}g,\end{equation}
\begin{equation} \label{rel19}
x_g\cdot y_{a,b}=x_{g^{[a,b]}}\cdot y_{a^{g^{[a,b]}},b^{g^{[a,b]}}}
\end{equation} for any $x_{g}\in
\mathcal X_O$ and any $y_{a,b}\in \mathcal Y_G$.\end{prop}
 \proof Follows from Claims
\ref{change}, \ref{fixed} and Lemmas
\ref{import}, \ref{imp}. \qed \\

Since every admissible automorphism group $H_{n,p}$ contains the
group $Br_{n,p}$, for any collection $\mathcal H$ of admissible
automorphism groups there is a natural epimorphism $r_{\mathcal
H-\text{equiv}}:\mathbb S(G,O)\to \mathbb S_{\mathcal
H-\text{equiv}}(G,O)$ of semigroups.

The semigroup $\mathbb S(\mathcal S_d,\mathcal S_d\setminus \{
{\bf{1}}\})$ that we denote by $V\mathbb S_d$ will be called {\it a
strong} {\it algebraic versal degree $d$ covering semigroup}. Note
that an embedding $i:G\hookrightarrow \mathcal S_d$ of a group $G$
into $\mathcal S_d$ induces the semigroup embedding of $\mathbb
S(G,O)$ into  $V\mathbb S_d$.

\begin{claim} The map $\alpha : Z\to G$ given by $\alpha (x_g)=g$
for $x_g\in X_O$ and $\alpha(y_{a,b})=[a,b]$ for $y_{a,b}\in Y_G$
defines a homomorphism $\mathbb S_{\mathcal H-\text{equiv}}(G,O)\to
G$.
\end{claim}
\proof Straightforward inspection of relations (\ref{rel12'}) and
(\ref{admis})  shows that for each of these relations the product of
the images of the left-side factors is equal in $G$ to the product
of the right-side factors.
\qed \\

Further on we denote this homomorphism $\mathbb S_{\mathcal
H-\text{equiv}}(G,O)\to G$ by  $\alpha_{G,\mathcal H-\text{equiv}}$,
or simply $\alpha_G$, and call it {\it the product homomorphism}.

The action $\rho$ of the group $G$ on the set $Z$, given by
$$ x_{g_1}\in X_O\mapsto \rho (g)(x_{g_1})=x_{gg_1g^{-1}}\in X_O,$$
$$y_{a,b}\in Y_G\mapsto \rho (g)(y_{a,b})=y_{gag^{-1},gbg^{-1}}\in Y_G,$$
 defines a homomorphism $\rho_{\mathbb S} : G\to
\text{Aut} (\mathbb S(G,O))$ and homomorphisms $\rho_{\mathcal
H-\text{equiv}} : G\to \text{Aut} (\mathbb S_{\mathcal
H-\text{equiv}}(G,O))$. Obviously, these actions are compatible with
the homomorphism $r_{\mathcal H-\text{equiv}}$.

If it does not lead to a confusion, we replace the notation $\mathbb
S_{\mathcal H-\text{equiv}}(G,O)$ by $\overline{\mathbb S}(G,O)$ and
then denote the both homomorphisms $\rho_{\mathbb S}$ and $\rho_{\mathcal
H-\text{equiv}}$ simply by $\rho$. The action $ \rho (g)$ on
$\overline{\mathbb S}(G,O)$ is called the {\it simultaneous
conjugation } by $g\in G$. Put $\lambda (g)=\rho(g^{-1})$ and
$\lambda _{\overline {\mathbb S}} =\lambda \circ\alpha_G$, $\rho _{\overline
{\mathbb S}} =\rho \circ\alpha_G$.

Whatever is an admissible  covering semigroup $\mathbb S_{\mathcal
H-\text{equiv}}(G,O)= \overline{\mathbb S}(G,O)$, the collection
$(\overline{\mathbb S}(G,O),G,\alpha_G, \rho)$ is a semigroup over
the group $G$ and the embedding $i:X_O\hookrightarrow Z$ defines an
embedding $i_*:S(G,O)\hookrightarrow \overline{\mathbb S}(G,O)$,
which is a semigroup homomorphism over $G$. Note also that
epimorphisms $r_{\mathcal H-\text{equiv}}:\mathbb S(G,O) \to \mathbb
S_{\mathcal H-\text{equiv}}(G,O)$ are also semigroup homomorphisms
over $G$.

Using relations (\ref{rel11}) -- (\ref{rel19}), any element $s\in
\overline{\mathbb S}(G,O)$ can be written in a so called {\it
reduced} form, $s=s_1\cdot s_2$, where $s_1\in S(G,O)$ and
$s_2=y_{a_1,b_1}\cdot \, .\, .\, .\, \cdot y_{a_p,b_p}$ for some
$a_1,b_1,\dots, a_p,b_p\in G$. We put  $\tau(s)=\tau(s_1)$ and
$g(s)=p$ and call them {\it type of $s$} and {\it genus of $s$},
respectively. It is easy to see that the type and the genus of $s\in
\overline{\mathbb S}(G,O)$ are well defined, that is, $\tau(s)$ and
$g(s)$ do not depend on the reduction of $s$ to a reduced form
$s=s_1\cdot s_2$.

Let $s_1\cdot s_2$ with $s_1=x_{g_1}\cdot \, .\, .\, .\, \cdot
x_{g_n}$ and  $s_2=y_{a,b}\cdot\, .\, .\, .\, \cdot y_{a_k,b_k}$ be
a reduced form of an element $s=s_1\cdot s_2 \in\overline{\mathbb
S}(G,O)$.  As it follows from Claim \ref{ggg}, the subgroup of $G$
generated by $g_1,\dots, g_n, a_1,b_1, \dots, a_k,b_k$ does not
depend on the choice of a reduced form of $s$. In what follows, we
denote this subgroup by $G_s$.

As in the case of factorization semigroups, for subgroups $H_1$ and
$H_2$ of a group $G$, we put
$$\overline{\mathbb S}(G,O)^{H_1}=\{s\in \overline{\mathbb S}(G,O)\mid G_s=H_1\},$$
$$\overline{\mathbb S}(G,O)_{H_2}= \{ s\in \overline{\mathbb S}(G,O)\mid \alpha (s)\in H_2 \},$$ and
$\overline{\mathbb S}(G,O)_{H_1}^{H_2}= \overline{\mathbb
S}(G,O)_{H_2}\cap \overline{\mathbb S}(G,O)^{H_1}$.

Let $G_{\Gamma}=(\widetilde G,\widetilde O)$ be the $C$-group
equivalent to $(G,O)$ (see subsection \ref{sub3}). For any set
$\mathcal R_H$ of admissible relations, the epimorphism
$\beta=\beta_{(G,O)}: (\widetilde G,\widetilde O)\to (G,O)$ of
equipped groups defines an epimorphism $\beta_*=\beta_{(G,O)*}:
\overline{\mathbb S}(\widetilde G,\widetilde O)=\mathbb
S_{H-\text{equiv}}(\widetilde G,\widetilde O)\to \overline{\mathbb
S}(G,O)=\mathbb S_{H-\text{equiv}}(G,O)$ over groups.

\begin{claim} \label{bb} The restriction of $\beta_*$ to the subsemigroup $S(\widetilde
G,\widetilde O)\subset \overline{\mathbb S}(\widetilde G,\widetilde
O)$ coincides with the isomorphism  of semigroups $S(\widetilde
G,\widetilde O)$ and $S(G,O)\subset \overline{\mathbb S}(G,O)$ {\rm
(}defined in subsection {\rm \ref{sub3})}.\end{claim} \proof
Obvious. \qed

\subsection{Solvability of some equations in strong covering semigroups}
\label{subs5} In this subsection, we will assume that in $\mathbb
S(G,O)$ there is the unity, ${\bf{1}}\in \mathbb S(G,O)$ (we add it
into $\mathbb S(G,O)$).

Let $s_1,s_2,s_3\in \mathbb S(G,O)$. We say that an equation
\begin{equation} \label{equa} s_1=s_2\cdot z\cdot s_3\end{equation}
is {\it solvable} in $S(G,O)\subset \mathbb S(G,O)$ if there is an
element $s\in S(G,O)$ such that $s_1=s_2\cdot s\cdot s_3$.

Note that $s$ is a solution of equation (\ref{equa}) if and only if
the following holds: if we write $s,s_1,s_2,s_3$ as products of
generators of $\mathbb S(G,O)$, then there is a finite sequence of
elementary transformations transforming the factorization of $s_1$
into the factorization of $s_2\cdot s\cdot s_3$, here an {\it
elementary transformation} means a change of some pair of two
neighboring factors into another one according to the one of
relations (\ref{rel11}) -- (\ref{rel19}) (reading either from the
left to the right or from the right to the left).

Consider four elements $s_1,\dots, s_4$ of  $\mathbb S(G,O)$ and let
us fix their presentations as products of generators of $\mathbb
S(G,O)$. Let  $S$ be a subset of  $S(G,O)$ the elements of which
have a fixed type. We say that the equations
\begin{equation} \label{equa1} s_{1}\cdot s\cdot s_2=s_{3}\cdot z\cdot
s_{4},\end{equation} where  $s\in S$, are {\it universally solvable}
if, first, there is a solution $\overline s\in S$ of (\ref{equa1})
for any $s$ and, second,  there is a finite sequence of elementary
transformations which satisfy the following property: for any
presentation of $s$ as the product of generators of $S(G,O)$ there
is a presentation of a solution $\overline s$ as the product of
generators such that this sequence of elementary transformations
transforms the factorization of $s_{1}\cdot s\cdot s_2$ into the one
of $s_{3}\cdot \overline s\cdot s_4$. The element $\overline s$
together with its factorization mentioned above will be called the
{\it universal solution} (for $s$ and its factorization) of equation
(\ref{equa1}).

\begin{claim} \label{cl0} For any $s\in S(G,O)$ of any fixed type and any $a,b\in G$ each of the following
equations
\begin{equation} \label{www1} s\cdot y_{a,b}=z\cdot y_{\alpha(s)a,b},\end{equation}
\begin{equation} \label{www2}  y_{a,b}\cdot s= y_{a,(\alpha(s))^{-1}b}\cdot z,\end{equation}
\begin{equation} \label{www3} s\cdot y_{a,b}=z\cdot y_{a,(\alpha(s)^{-1})^{[a,b]}b}\end{equation}
is universally solvable in $S(G,O).$\end{claim}

\proof Let $s=x_{g_1}\cdot \, .\, .\, .\, \cdot x_{g_n}$. If $n=1$,
then a universal solvability of equation (\ref{www1}) follows from
relation (\ref{rel13}). Assume that equation (\ref{www1}) is
 universally solvable for any $s$ of length $n-1$ and let us write an element $s$
of length $n$ in the form: $s=s_1\cdot x_{g_n}$. We have
$$s\cdot y_{a,b}=s_1\cdot x_{g_n}\cdot y_{a,b}=s_1\cdot z_1\cdot y_{g_na,b}=
\rho_{\overline{\mathbb S}}(s_1)(z_1)\cdot s_1\cdot y_{g_na,b},$$
where $z_1$ is an  universal solution of equation $x_{g_n}\cdot
y_{a,b}=z\cdot y_{g_na,b}$. By assumption, for some $z_2$ we have:
$s_1\cdot y_{g_na,b}=z_2\cdot y_{\alpha(s_1)ga,b}$. Therefore
$$s\cdot y_{a,b}=\rho_{\overline{\mathbb S}}(s_1)(z_1)\cdot s_1\cdot
y_{g_na,b}=\rho_{\overline{\mathbb S}}(s_1)(z_1)\cdot z_2\cdot
y_{\alpha(s_1)g_na,b}=(\rho_{\overline{\mathbb S}}(s_1)(z_1)\cdot
z_2)\cdot y_{\alpha(s)a,b},
$$ that is,  equation (\ref{www1}) is universally solvable always.

The proof of universal solvability of equation (\ref{www2}) is
similar to one for equation (\ref{www1}). Only we must use relation
(\ref{rel16}) instead of relation (\ref{rel13}).

To prove the universal solvability of equation (\ref{www3}), note
that $\alpha(\lambda([a,b])(s))=(\alpha(s))^{[a,b]}$. Therefore, by
relation (\ref{rel12}), we have
$$s\cdot y_{a,b}=y_{a,b}\cdot\lambda([a,b])(s)=
y_{a,(\alpha(s)^{-1})^{[a,b]}b}\cdot z_1,$$ where $z_1$ is a
universal solution of equation $y_{a,b}\cdot\lambda([a,b])(s)=
y_{a,(\alpha(s)^{-1})^{[a,b]}b}\cdot z$. Now to prove the universal
solvability of equation (\ref{www3}), it suffices several times to
use relation (\ref{rel12}). \qed

The following proposition is a generalization of Main Lemma 2.1 in
\cite{Ka}.

\begin{prop} \label{Kan} Let $O\subset G$ be a finite set. Then for any
$$h\in \langle g_1,\dots, g_n, a_1,b_1,\dots, a_k,b_k\rangle$$ and
for any $s_{g^{-1}}\in S(G,O)$ such that $\alpha(s_{g^{-1}})=g^{-1}$
the following equation
$$x_{g_1}\cdot .\, .\, .\cdot x_{g_{n}}\cdot x_g\cdot s_{g^{-1}}\cdot y_{a_1,b_1}\cdot .\, .\, .\cdot y_{a_k,b_k}=
x_{g_1}\cdot .\, .\, .\cdot x_{g_{n}}\cdot x_{g^h}\cdot z\cdot
y_{a_1,b_1}\cdot .\, .\, .\cdot y_{a_k,b_k}
$$ is solvable in
$S(G,O)$.
\end{prop}

\proof For $s_1, s_2\in \mathbb S(G,O)$ denote by $G(s_1,s_2)$ the
subset of $G$ such that $h\in G(s_1,s_2)$ if and only if the
equations \begin{equation} \label{eqxx} s_1\cdot x_g\cdot
s_{g^{-1}}\cdot s_2=s_1\cdot x_{g^h}\cdot z\cdot s_2\end{equation}
are universally solvable in $S(G,O)$ for each $g\in O\subset G$ and
any $s_{g^{-1}}\in S(G,O)$ of fixed type and such that
$\alpha(s_{g^{-1}})=g^{-1}$. Note that if $z_1$ is a solution of
equation (\ref{eqxx}), then $\alpha(z_1)=(g^h)^{-1}$.

\begin{claim} \label{sol1} Let $O\subset G$ be a finite set. Then for any $s_1, s_2\in
\mathbb S(G,O)$ the set $G(s_1,s_2)$ is a subgroup of $G$.
\end{claim}
\proof Let us show that if $h\in G(s_1,s_2)$, then $h^{-1}\in
G(s_1,s_2)$.

Let $z_1$ be an universal solution of equation (\ref{eqxx}). If we
apply the inverse sequence of the sequence of elementary
transformations giving the universe resolution of equation
(\ref{eqxx}), then it is easy to see that $s_{g^{-1}}$ is the
 universal solution of the equation
$$s_1\cdot x_{g_1}\cdot z_1\cdot s_2=s_1\cdot x_{g_1^{h_1}}\cdot
z\cdot s_2,$$ where $g_1=g^h$. For any $h\in G$ the conjugation of
the elements of $S(G,O)$ by $h$ is a bijection, and for each $g\in
G$ the set of different factorizations of elements $s\in S(G,O)$ of
fixed type and such that $\alpha(s)= g^{-1}$ is finite. Therefore
$h^{-1}\in G(s_1,s_2)$ if $h\in G(s_1,s_2)$.

Let us show that if $h_1,h_2\in G(s_1,s_2)$, then $h_1h_2\in
G(s_1,s_2)$. It is easy to see that if $z_1$ is an universal
solution of the equation $$s_1\cdot x_g\cdot s_{g^{-1}}\cdot
s_2=s_1\cdot x_{g^{h_1}}\cdot z\cdot s_2$$ and if $z_2$ is an
universal solution of the equation $$s_1\cdot x_{g^{h_1}}\cdot
z_1\cdot s_2=s_1\cdot x_{(g^{h_1})^{h_2}}\cdot z\cdot s_2,$$ then
$z_2$ is an universal solution of the equation
$$s_1\cdot x_{g}\cdot s_{g^{-1}}\cdot s_2=s_1\cdot x_{g^{h_1h_2}}\cdot z\cdot
s_2. \hspace{3cm} \qed $$

\begin{claim} \label{sol2} For any $s_1=s_1^{\prime}\cdot s_1^{\prime\prime}, s_2=s_2^{\prime}\cdot s_2^{\prime\prime}\in
\mathbb S(G,O)$, we have $G(s_1^{\prime\prime},s_2^{\prime})\subset
G(s_1,s_2)$.
\end{claim}
\proof Obvious. \qed

\begin{claim} \label{sol3} For  any $h_1, \dots , h_n\in O$, we have
$\langle h_1,\dots, h_n\rangle\subset G(x_{h_1}\cdot .\, .\, .\, .
\cdot x_{h_n},\bf{1})$.
\end{claim}
\proof It easily follows from Claim \ref{sol1} and from the
following equalities: $$x_h\cdot x_g\cdot
s_{g^{-1}}=x_{g^{(h^{-1})}}\cdot \rho(h)(s_{g^{-1}})\cdot
x_h=x_h\cdot x_{g^{(h^{-1})}}\cdot \rho(h)(s_{g^{-1}}), $$ since
 $x_{g^{(h^{-1})}}\cdot
\rho(h)(s_{g^{-1}})\cdot x_h=\rho_{\overline{\mathbb
S}}(x_{g^{(h^{-1})}}\cdot \rho(h)(s_{g^{-1}}))(x_h)\cdot
x_{g^{(h^{-1})}}\cdot \rho(h)(s_{g^{-1}})$ and
$\alpha(x_{g^{(h^{-1})}}\cdot \rho(h)(s_{g^{-1}}))={\bf{1}}$. \qed
\begin{claim} \label{sol4} For  any $a,b\in G$, we have $ab^{-1}a^{-1}\in G({\bf{1}},y_{a,b})$.
\end{claim}
\proof  By Claim \ref{cl0}, we have $$x_g\cdot s_{g^{-1}}\cdot
y_{a,b}=x_g\cdot z_1\cdot y_{g^{-1}a,b}=\rho(g)(z_1)\cdot x_g\cdot
y_{g^{-1}a,b},$$ where $z_1$ is an universal solution of equation
$s_{g^{-1}}\cdot y_{a,b}=z\cdot y_{g^{-1}a,b}$. By relation
(\ref{rel13}), $$x_g\cdot
y_{g^{-1}a,b}=x_{gg^{-1}aba^{-1}ggg^{-1}ab^{-1}a^{-1}}\cdot
y_{a,b}=x_{g^{ab^{-1}a^{-1}}}\cdot y_{a,b}.$$  Therefore
$$x_g\cdot
s_{g^{-1}}\cdot y_{a,b}=\rho(g)(z_1)\cdot x_g\cdot
y_{g^{-1}a,b}=\rho(g)(z_1)\cdot x_{g^{ab^{-1}a^{-1}}}\cdot
y_{a,b},$$ that is, $ab^{-1}a^{-1}\in G({\bf{1}},y_{a,b})$. \qed

\begin{claim} \label{sol5} For  any $a,b\in G$, we have $aba^{-1}b^{-1}a^{-1}\in G({\bf{1}},y_{a,b})$.
\end{claim}
\proof  Applying $ln(x_g\cdot s_{g^{-1}})$ times relation
(\ref{rel12}) and after that applying $ln(s_{g^{-1}})$ times
relation (\ref{rel11}), we have
$$x_g\cdot s_{g^{-1}}\cdot y_{a,b}=
y_{a,b}\cdot x_{g^{[a,b]}}\cdot z_1=y_{a,b}\cdot z_2\cdot
x_{g^{[a,b]}},$$ where $z_1=\lambda([a,b])(s_{g^{-1}})$ and
$z_2=\rho(g^{[a,b]})(z_1)$. It is easy to see that
$\alpha(z_2)=(g^{[a,b]})^{-1}$. By Claim \ref{cl0} and relation
(\ref{rel16}),
$$y_{a,b}\cdot z_2\cdot
x_{g^{[a,b]}}=z_3\cdot y_{a,g^{[a,b]}b}\cdot x_{g^{[a,b]}}=z_3\cdot
y_{a,b}\cdot x_{g^{[a,b]ba^{-1}b^{-1}}},$$ where $z_3$ is an
universal solution of equation $y_{a,b}\cdot z_2=z\cdot
y_{a,g^{[a,b]}b}$.

Applying relation (\ref{rel12}), we obtain $$z_3\cdot y_{a,b}\cdot
x_{g^{[a,b]ba^{-1}b^{-1}}}=z_3\cdot
x_{g^{[a,b]ba^{-1}b^{-1}[b,a]}}\cdot y_{a,b}=z_3\cdot
x_{g^{aba^{-1}b^{-1}a^{-1}}}\cdot y_{a,b}$$ and to complete the
proof it suffices to use the relation $$z_3\cdot
x_{g^{aba^{-1}b^{-1}a^{-1}}}=x_{g^{aba^{-1}b^{-1}a^{-1}}}\cdot
\lambda(g^{aba^{-1}b^{-1}a^{-1}})(z_3). \qed$$

\begin{claim} \label{sol6} For  any $a,b\in G$, we have $a,b\in G({\bf{1}},y_{a,b})$.
\end{claim}
\proof By Claims \ref{sol4} and \ref{sol5}, the elements
$ab^{-1}a^{-1}$ and $aba^{-1}b^{-1}a^{-1}$ belong to
$G({\bf{1}},y_{a,b})$. It follows from Claim \ref{sol1} that
$$(ab^{-1}a^{-1})(aba^{-1}b^{-1}a^{-1})=(ab)^{-1}\in
G({\bf{1}},y_{a,b}).$$ Therefore $ab\in G({\bf{1}},y_{a,b})$ and
hence $a=(ab^{-1}a^{-1})(ab)\in G({\bf{1}},y_{a,b})$. Applying one
more Claim \ref{sol1} to $a$ and $ab$, we obtain that $b\in
G({\bf{1}},y_{a,b})$. \qed

\begin{claim} \label{sol7} For  any $a_1,b_1, \dots, a_k,b_k\in G$,
we have $$\langle a_1,b_1,\dots, a_k,b_k\rangle\subset
G({\bf{1}},y_{a_1,b_1}\cdot \, .\, .\, . \cdot y_{a_k,b_k}).$$
\end{claim}

\proof Let us use induction on $k$. In the case $k=1$, it is Claim
\ref{sol6}. Assume that for some $k-1$ Claim is true. Therefore, by
Claim \ref{sol2}, $$\langle a_1,b_1,\dots,
a_{k-1},b_{k-1}\rangle\subset G({\bf{1}},y_{a_1,b_1}\cdot \, .\, .\,
. \cdot y_{a_{k},b_{k}}).$$ Denote by $u_t=[a_1,b_1]\dots
[a_t,b_t]$. We have
$$x_g\cdot s_{g^{-1}}\cdot y_{a_1,b_1}\cdot \, .\, .\, . \cdot y_{a_{k-1},b_{k-1}}\cdot y_{a_k,b_k}
=y_{a_1,b_1}\cdot \, .\, .\, . \cdot y_{a_{k-1},b_{k-1}}\cdot
x_{g^{u_{k-1}}}\cdot \lambda(u_{k-1})(s_{g^{-1}})\cdot
y_{a_k,b_k}.$$ Denote by $z_1=\lambda(u_{k-1})(s_{g^{-1}})$ and let
$z_c$ (see Claim \ref{sol6}) be an universal solution of equation
$x_{g^{u_{k-1}}}\cdot z_1\cdot y_{a_k,b_k} =x_{g^{u_{k-1}c}}\cdot
z_c\cdot y_{a_k,b_k}$, where $c= a_k$ or $b_k$. We have
$$
\begin{array}{c} x_g\cdot s_{g^{-1}}\cdot y_{a_1,b_1}\cdot \, .\, .\, . \cdot
y_{a_{k},b_{k}}= \\
y_{a_1,b_1}\cdot \, .\, .\, . \cdot y_{a_{k-1},b_{k-1}}\cdot
x_{g^{u_{k-1}c}}\cdot  z_c\cdot y_{a_k,b_k}=
\\
  x_{g^{c^{u_{k-1}}}}\cdot \lambda(u_{k-1})(z_c)\cdot
y_{a_1,b_1}\cdot \, .\, .\, . \cdot y_{a_{k-1},b_{k-1}}\cdot
y_{a_k,b_k}.  
\end{array}
$$
Now, since by assumption,  $u_{k-1}\in G({\bf{1}},y_{a_1,b_1}\cdot
\, .\, .\, . \cdot y_{a_k,b_k})$, we obtain that the element $c\in
G({\bf{1}},y_{a_1,b_1}\cdot \, .\, .\, . \cdot y_{a_k,b_k})$. \qed

Now the proof of Proposition \ref{Kan} easily follows from Claims
\ref{sol1} -- \ref{sol7}. \qed

\begin{prop} \label{prop-main}
Let $(G,O)$ be an equipped finite group such that the elements of
$O=C_1\sqcup\dots \sqcup C_m$ generate the group $G$. Denote by
$n_i$ the number of elements of the conjugacy class $C_i$ and by
$p_i$ the order of elements of $C_i$. Let $s_1\in S(G,O)$ be such
that $\tau_i(s_1)> n_ip_i$ for all $i$, $1\leq i\leq m$, and
$s_2=y_{a_1,b_1}\cdot .\, .\, .\cdot y_{a_k,b_k}$ be such that
$s_1\cdot s_2\in \mathbb S(G,O)^G$. Then the equation
$$s_1\cdot s_2=z\cdot (y_{{\bf 1},{\bf{1}}})^k$$ is solvable in $S(G,O)$. \end{prop}

\proof  Let $a_1=g_1^{-1}\dots g_n^{-1}$, where $g_1,\dots,g_n\in
O$, and let $s_1=x_{h_1}\cdot .\, .\, .\, . \cdot x_{h_N}$. Let
$g_1\in C_i$. Since $\tau_i(s_1)> n_ip_i$, among the factors
$x_{h_1},\dots,x_{h_N}$ there are at least $p_i+1$ factors with the
same $h_j\in C_i$. Moving $p_i$ of these factors to the right (using
relation (\ref{rel11})), we obtain that
$$s_1\cdot s_2=s_1'\cdot x_{h_j}\cdot (x_{h_j})^{p_i-1}\cdot s_2$$
for some $s_1'\in S(G,O)$ such that $s'_1\cdot s_2\in
\overline{\mathbb S}(G,O)^G$.

Applying Proposition \ref{Kan}, we have
$$s_1'\cdot x_{h_j}\cdot (x_{h_j})^{p_i-1}\cdot s_2=s_1'\cdot s'\cdot x_{g_1}\cdot s_2$$
for some $s'\in S(G,O)$ such that $\tau(s')=\tau((x_{h_j})^{p_i-1})$
and $\alpha(s')=g_1^{-1}$.

By relation (\ref{rel13}), we have
$$s_1'\cdot s'\cdot x_{g_1}\cdot s_2=s_1'\cdot s'\cdot x_{g_1}\cdot y_{a_1,b_1}\cdot .\, .\, .\cdot
y_{a_k,b_k}=s_1'\cdot s'\cdot x_{g'_1}\cdot y_{a'_1,b_1}\cdot .\,
.\, .\cdot y_{a_k,b_k},$$ where $g'_1$ is an element conjugate to
$g_1$ and $a'_1=g_2^{-1}\dots g_n^{-1}$.

Note that $\widetilde s_1=s_1'\cdot s'\cdot x_{g'_1}$ and
$\widetilde s_2=y_{a'_1,b_1}\cdot .\, .\, .\cdot y_{a_k,b_k}$
satisfy all conditions of Proposition \ref{prop-main}. Therefore, by
induction on $n$, we obtain that $$s_1\cdot s_2=\overline s_1\cdot
(y_{{\bf 1},b_1}\cdot y_{a_2,b_2}\cdot .\, .\, .\cdot y_{a_k,b_k})$$
for some $\overline s_1$ which together with $\overline s_2=y_{{\bf
1},b_1}\cdot y_{a_2,b_2}\cdot .\, .\, .\cdot y_{a_k,b_k}$ satisfies
all conditions of Proposition \ref{prop-main}.

The same arguments (only instead of relation (\ref{rel13}) we must
use relation (\ref{rel16})) give that $$s_1\cdot s_2=\widetilde
s_1\cdot (y_{{\bf 1},{\bf 1}}\cdot y_{a_2,b_2}\cdot .\, .\, .\cdot
y_{a_k,b_k})$$ for some $\widetilde s_1$ which together with
$\widetilde s_2=y_{a_2,b_2}\cdot .\, .\, .\cdot y_{a_k,b_k}$
satisfies all conditions of Proposition \ref{prop-main}.

Now to complete the proof of Proposition \ref{prop-main}, it
suffices to note that $y_{{\bf 1},{\bf 1}}$ commutes with any
element $s\in S(G,O)$ (relation (\ref{rel12})) and to use induction
on $k$. \qed

\subsection{ On the number of solutions of equation $\alpha(z)=g$}
\label{number} Let $(G,O)$ be an equipped group such that the
elements of $O$ generate the group $G$, $O=C_1\sqcup \dots \sqcup
C_m$ decomposition into a disjoint union of conjugacy classes of
$G$. In this subsection, we investigate the following problem: {\it
for fixed type $\overline t\in \mathbb Z_{\geq 0}^m$, fixed genus
$p$ and given element $h\in G$ to estimate the number of solutions
in an admissible covering semigroup $\overline{\mathbb S}(G,O)^G$ of
the equation $\alpha_G(z)=h$ under the
restrictions
$\tau(z)=\overline t$ and $g(z)=p$. }

In \cite{Ku3}, this problem was solved in the case of $p=0$ and
$\overline t=(t_1,\dots, t_m)$ such that all $t_i$ are big enough
(see Theorems \ref{Cuniq} and \ref{Cuniq2}). In this subsection we
generalize these  results  to the case of arbitrary genus, namely,
we prove

\begin{thm}  \label{Cuniq3} Let $\overline{\mathbb S}(G,O)=\mathbb S_{H-\text{equiv}}(G,O)$
be an admissible covering semigroup over
an equipped finite group $(G,O)$, $O=C_1\sqcup \dots \sqcup C_m$.
Then there is a constant $T_1\in \mathbb N$ such that if  for an
element $\overline s_1\in \overline{\mathbb S}(G,O)^G$ the $i$-th
type $\tau_i(s_1)\geq T_1$ for all $i=1,\dots, m$, then there are
$a_{(G,O)}$ elements $\overline s_1,\dots, \overline
s_{a_{(G,O)}}\in \overline{\mathbb S}(G,O)^G$ such that
\begin{itemize}
\item[$(1)$] $\overline s_i\neq \overline s_j$ for $1\leq i< j\leq a_{(G,O)}$;
\item[$(2)$] $\tau(\overline s_i)=\tau(\overline s_1)$ for $1\leq i\leq a_{(G,O)}$;
\item[$(3)$]$g(\overline s_i)=g(\overline s_1)$ for $1\leq i\leq a_{(G,O)}$;
\item[$(4)$] $\alpha_G(\overline s_i)=\alpha_G(\overline s_1)$ for $1\leq i\leq a_{(G,O)}$;
\item[$(5)$] if $\overline s\in \overline S(G,O)^G$ is such that $\tau(\overline s)=\tau(\overline s_1)$,
$g(\overline s)=g(\overline s_1)$, and $\alpha_G(\overline
s)=\alpha_G(\overline s_1)$, then $\overline s=\overline s_i$ for
some $i$, $1\leq i\leq a_{(G,O)}$.
\end{itemize}
\end{thm}
\proof Let $p$ be the genus of $\overline s_1$ and $T$ be a constant
the existence of which is claimed in Theorem \ref{Cuniq}. Without
loss of generality, we can assume that $T> \max_{1\leq i\leq
m}n_ip_i$, where $n_i$ is the number of elements of the conjugacy
class $C_i$ and $p_i$ is the order of elements of $C_i$. By
Proposition \ref{prop-main} (and since
$r_{\mathcal H-\text{equiv}}$
is an epimorphism),
the element $\overline s_1$ can be written in the form $\overline
s_1=s_1\cdot (y_{{\bf 1},{\bf 1}})^{p}$. By Theorem \ref{Cuniq},
there are exactly $a_{(G,O)}$ different elements $s_1,\dots,
s_{a_{(G,O)}}\in S(G,O)^G$ satisfying conditions (1) -- (4) of
Theorem \ref{Cuniq3}.

Consider the elements $\overline s_i=s_i\cdot (y_{{\bf 1},{\bf
1}})^{p}$, $1\leq i\leq a_{(G,O)}$. By Proposition \ref{prop-main}
and  Theorem \ref{Cuniq}, they satisfy conditions (2) -- (5) of
Theorem \ref{Cuniq3}. Let us show that they also satisfy condition
(1)  of Theorem \ref{Cuniq3}. Assume that for some $i\neq j$ we have
$\overline s_i=\overline s_j$, that is, if we write $s_i$ and $s_j$
as products of generators $x_g$, $g\in O$, then there is a finite
sequens  $\text{Tr}$ of elementary transformations  (see subsection
\ref{subs5}) transforming the factorization of $\overline s_i$ into
the factorization of $\overline s_j$. By Claim \ref{bb}, the
selected factorizations alow us to lift the element $\overline s_i$
into $\overline{\mathbb S}(\widetilde G,\widetilde O)=\mathbb
S_{H-\text{equiv}}(\widetilde G,\widetilde O)$ (that is, to consider
an element $\widetilde s_i\in \overline{\mathbb S}(\widetilde
G,\widetilde O)$ with the same factorization as the one of
$\overline s_i$), where $G_{\Gamma}=(\widetilde G,\widetilde O)$ is
the $C$-group equivalent to $(G,O)$. Let us apply the same sequence
$\text{Tr}$ of elementary transformations to the element $\widetilde
s_i$. As a result we obtain an element $\widetilde s_j=s_j\cdot
(y_{a_1,b_1}\cdot .\, .\, .\cdot y_{a_p,b_p})=\widetilde s_i$ such
that $a_l, b_l\in \ker \beta_{(G,O)}\subset Z(\widetilde G)$ for
$1\leq l\leq p$. But it is impossible, since in this case we have
$\alpha_{\widetilde G}(y_{a_l,b_l})={\bf{1}}$ for $1\leq l\leq p$
and therefore by Remark \ref{rem}, we must have $\alpha_{\widetilde
G}(\widetilde s_i)\neq \alpha_{\widetilde G}(\widetilde s_j)$. \qed

\begin{thm}  \label{Cuniq4} Let $G$ be a finite group and
$O^{\prime}\subset O$ be two its equipments such that the elements
of $O^{\prime}=C_1\sqcup \dots \sqcup C_k$ generate the group $G$.
Then there is a constant $T=T_{O^{\prime}}$ such that if for an
element $\overline s_1\in \overline{\mathbb S}(G,O)^G$ its $i$-th
type $\tau_i(\overline s_1)\geq T$ for all $i=1,\dots, k$, then
there are not more than $a_{(G,O^{\prime})}$ elements $\overline
s_1,\dots, \overline s_{n}\in \overline{\mathbb S}(G,O)^G$ such that
for $1\leq i< j\leq n$
\begin{itemize}
\item[$(i)$] $\overline s_i\neq \overline s_j$;
\item[$(ii)$] $\tau(\overline s_i)=\tau(\overline s_1)$;
\item[$(iii)$] $\alpha_G(\overline s_i)=\alpha_G(\overline s_1)$,
\item[$(iv)$]$g(\overline s_i)=g(\overline s_1)$,
\end{itemize}
where $a_{(G,O^{\prime})}$ is the ambiguity index of
$(G,O^{\prime})$.
\end{thm}
\proof It is similar to the proof of Theorem \ref{Cuniq3}.\qed

\begin{cor} \label{cor1} Let $C_1$ be a conjugacy class of an odd permutation
$\sigma_1\in \mathcal S_d$ such that $\sigma_1$ leaves fixed at
least two elements. Then in the case when $C_1$ is contained in an
equipment $O=C_1\sqcup\dots\sqcup C_m$ of $\mathcal S_d$, there is a
constant $T=T_{C_1}$ such that for any $\sigma\in \mathcal S_d$, any
fixed integer $p\geq 0$, and any $\overline t=(t_1,\dots, t_m)\in
\mathbb Z^m_{\geq 0}$ such that $t_1\geq T$ the equation
\begin{equation} \label{eqxxx} \alpha_{\mathcal S_d}(z)=\sigma\end{equation}  has in each covering
semigroup $\overline{\mathbb S}(\mathcal S_d,O)$ at most one
solution $\overline s$ satisfying conditions $g(\overline s)=p$ and
$\tau(\overline s)=\overline t$. Under assumption $t_1\geq T$, the
existence of solution of equation {\rm (\ref{eqxxx})} does not
depend on $p$ and depends only on $t$.\end{cor} \proof It follows
from Theorem
\ref{Cuniq4} and the main result of \cite {Ku2}. \qed\\

Let $\overline{\mathbb S}(G,O)$ be an admissible covering semigroup
over an equipped finite group $(G,O)$, $O=C_1\sqcup \dots \sqcup
C_m$, and let $T=(t_1,\dots, t_m)\in\mathbb Z_{\geq o}$. Denote by
$$\overline{\mathbb S}(G,O)_{\geq T}=\{ \overline s\in
\overline{\mathbb S}(G,O)\mid \tau_i(\overline s)\geq t_i\}$$ a
subsemigroup of $\overline{\mathbb S}(G,O)$. By Theorem
\ref{Cuniq3}, we have
\begin{cor}\label{cor2} For any equipped finite group $(G,O)$,
$O=C_1\sqcup \dots \sqcup C_m$, there is a constant $T=(t_1,\dots,
t_m)\in\mathbb Z_{\geq o}$ such that the restrictions
$r_{H-\text{equiv}}:\mathbb S(G,O)^G_{\geq T}\to \mathbb
S_{H-\text{equiv}}(G,O)^G_{\geq T}$ of the epimorphisms
$r_{H-\text{equiv}}:\mathbb S(G,O)\to \mathbb
S_{H-\text{equiv}}(G,O)$ to the subsemigroup $\mathbb S(G,O)^G_{\geq
T}$ are isomorphisms for any set $\mathcal R_{H}$ of admissible
relations.
\end{cor}

\section{Geometric semigroups of coverings} \label{secG}

\subsection{Monodromy encoding of ramified coverings}\label{encoding-ramified}
To describe ramified coverings of a given connected manifold $M$ (in
most cases, $M$ will be a connected compact oriented surface with
one hole and a base point on the boundary) we use  a traditional
{\it monodromy encoding} of non-ramified coverings and the unicity
of {\it manifold ramified completions}.

Namely, given a non-ramified, possibly disconnected, degree $d$
covering $\pi: \tilde M\to M$ over a connected manifold $M$ with a
base point $q\in M$, the lifts of a loop at $q$ form a set of $d$
paths in $\tilde M$ starting each at a different point of
$\pi^{-1}(q)$ and, thus, they give rise to a permutation of the set
$\pi^{-1}(q)$. This permutation depends only on the homotopy type of
the loop and in this way one obtains an encoding of the covering by
{\it a homomorphism from $\pi_1(M,q)$ to the permutation group of
$\pi^{-1}(q)$}.  In particular, if the covering $\pi$ is equipped
with a marking $\nu :  I_d=[1,\dots,d]\to \pi^{-1}(q)$ it gives a
well defined homomorphism $\pi_1(M,q)\to \mathcal S_d$.

The foliowing Proposition is well known and straightforward.

\begin{prop}\label{encoding} Two non-ramified marked coverings over
the same based space $(M,q)$ are isomorphic as marked coverings  if
and only if they define the same homomorphism $\pi_1(M,q)\to
\mathcal S_d$. If an isomorphism exists, it is unique; equivalently,
marked coverings have no non-trivial automorphisms.

Each homomorphism $\pi_1(M,q)\to \mathcal S_d$ corresponds to a
certain non-ramified marked degree $d$ covering of $M$. The covering
space is connected if and only if  the action of $\pi_1(M,q)$ on
$\pi^{-1}(q)$ is transitive. The orbits of the action correspond
canonically to connected components of the covering space.
\end{prop}

Manifold completions of non-ramified coverings by ramified ones are
most transparent in low dimensions.

The following result is also well known and straightforward; it can
be found, for example, in \cite{Fox}.

\begin{prop}\label{completion} Let $(M,q)$ be a based two-dimensional manifold and
$B=\{P_1,\dots, P_n\}$ a finite subset of $M$ disjoint from $q$ and
$\partial M$. Then each marked non-ramified covering of $(U,q)$ with
$U=M\setminus B$ has one and only one ramified completion $\tilde
M\to M$.
\end{prop}

\subsection{Surfaces with a hole and their skeletons}\label{skeletons} Our main
building blocks are connected compact oriented surfaces
($2$-dimensional manifolds) with one hole and one marked point on
the boundary. We equip them, in addition, with {\it a
semi-skeleton}, {\it skeleton}, or {\it caudate skeleton}.

We define a semi-skeleton of a connected compact oriented
$2$-dimensional manifold (with one hole) to be the union of disjoint
embedded bouquets of two oriented circles with a property that the
complement of the union is homeomorphic to a punctured disc.
Clearly, the number of connected components of a semi-skeleton is
equal to the genus of the surface. We distinguish the two circles of
a bouquet that represents a connected component of a semi-skeleton
by means of intersection index, namely, we speak on $\lambda$- and
$\mu$-{\it circles} in a bouquet $C_\lambda\vee C_\mu$ by respecting
the convention that $C_\lambda\cdot C_\mu=-1$ (and not $1$).

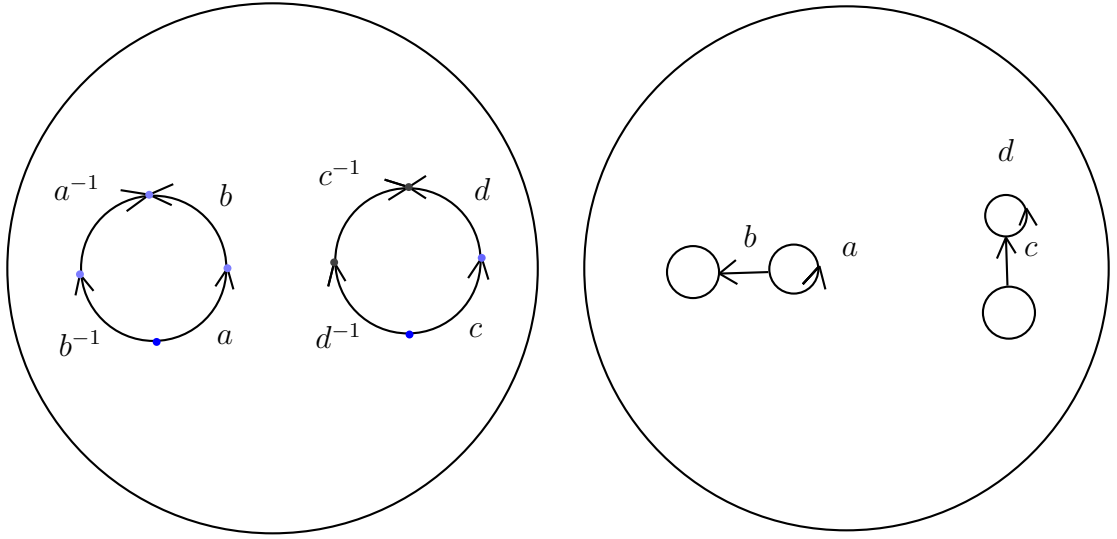
\begin{figure}\caption{Plane diagrams of a genus-$2$ surface with its semi-skeleton.}\label{glasses}
\vskip0.1in \pagestyle{empty}
\newrgbcolor{xdxdff}{0.49 0.49 1}
\newrgbcolor{wwwwff}{0.4 0.4 1}
\psset{xunit=1.0cm,yunit=1.0cm,algebraic=true,dotstyle=o,dotsize=3pt
0,linewidth=0.8pt,arrowsize=3pt 2,arrowinset=0.25}
\begin{pspicture*}(-3.55,-5.75)(25.01,7.7)
\pscircle(0.5,3.6){3.54} \pscircle(-1.08,3.6){0.98}
\pscircle(2.3,3.7){0.98} \psline(-2.06,3.52)(-1.88,3.3)
\psline(-2.06,3.52)(-2.12,3.24) \psline(-1.14,4.58)(-0.84,4.44)
\psline(-0.82,4.72)(-1.14,4.58) \psline(-0.1,3.6)(-0.26,3.36)
\psline(-0.02,3.32)(-0.1,3.6) \psline(-1.44,4.36)(-1.14,4.58)
\psline(-1.52,4.64)(-1.14,4.58) \rput[tl](-0.24,2.77){$a$}
\rput[tl](-0.21,4.73){$b$} \rput[tl](-2.41,4.84){$a^{-1}$}
\rput[tl](-2.34,2.77){$b^{-1}$} \psline(1.24,3.36)(1.32,3.68)
\psline(1.32,3.68)(1.24,3.36) \psline(1.47,3.43)(1.32,3.68)
\psline(2.04,4.51)(2.31,4.68) \psline(1.99,4.78)(2.31,4.68)
\psline(2.31,4.68)(1.99,4.78) \psline(2.31,4.68)(2.54,4.53)
\psline(2.31,4.68)(2.54,4.83) \psline(3.28,3.74)(3.37,3.45)
\psline(3.28,3.74)(3.1,3.47) \rput[tl](3.11,2.87){$c$}
\rput[tl](3.18,4.8){$d$} \rput[tl](1.11,5.04){$c^{-1}$}
\rput[tl](1.07,2.87){$d^{-1}$} \pscircle(8.13,3.6){3.5}
\pscircle(6.09,3.55){0.36} \pscircle(7.43,3.59){0.34}
\pscircle(10.25,4.3){0.29} \pscircle(10.29,3.01){0.36}
\psline(6.44,3.53)(7.09,3.55) \psline(10.25,4.01)(10.27,3.37)
\psline(6.44,3.53)(6.68,3.39) \psline(6.66,3.71)(6.44,3.53)
\psline(7.77,3.63)(7.55,3.39) \psline(7.77,3.63)(7.55,3.39)
\psline(7.77,3.63)(7.55,3.39) \psline(7.55,3.39)(7.77,3.63)
\psline(7.85,3.32)(7.77,3.63) \psline(10.13,3.73)(10.25,4.01)
\psline(10.38,3.8)(10.25,4.01) \rput[tl](8.07,3.94){$a$}
\rput[tl](6.76,4.18){$b$} \rput[tl](10.49,3.94){$c$}
\psline(10.66,4.17)(10.52,4.4) \psline(10.43,4.17)(10.52,4.4)
\rput[tl](10.14,5.32){$d$}
\begin{scriptsize}
\psdots[dotstyle=*,linecolor=blue](-1.04,2.62)
\psdots[dotstyle=*,linecolor=blue](2.32,2.72)
\psdots[dotstyle=*,linecolor=xdxdff](-0.1,3.6)
\psdots[dotstyle=*,linecolor=xdxdff](-1.14,4.58)
\psdots[dotstyle=*,linecolor=xdxdff](-2.06,3.52)
\psdots[dotstyle=*,linecolor=wwwwff](3.28,3.74)
\psdots[dotstyle=*,linecolor=darkgray](1.32,3.68)
\psdots[dotstyle=*,linecolor=darkgray](2.31,4.68)
\end{scriptsize}
\end{pspicture*}
\end{figure}

Given such a triple $(F,q, S^\infty)$, where $S^\infty$ is a
semi-skeleton of a surface $F$ with a fixed point $q\in\partial F$,
we can represent it by {\it an open-eyes plane diagram}, that is to
draw a disc with a marked point on its boundary and $p=g(F)$ holes
inside the disc, and trace the standard $4$-gone identification
scheme on the boundary of each hole, see the left drawing on Figure
\ref{glasses}; the orientation  of $F$ should be induced by the
counter-clock wise orientation of the disc. When it happens to be
more convenient and transparent, we use also another, equivalent,
presentation and draw a disc with "pince-nez", that is $p$ pairs of
holes with "bridges",  see the right drawing on Figure \ref{glasses}
(there, the $\lambda$-circles are $a$ and $c$, while $b$ and $d$ are
the $\mu$-circles).

Open-eyes plane diagrams of a given triple $(F,q,S^\infty)$ are
defined up to orientation preserving (stratified) homeomorphisms of
the disc respecting the marked point and the orientation of the
boundary identification strata. Converse statement is also true, an
open-eyes plane diagram defines the triple $(F,q,S^\infty)$ up to
orientation preserving homeomorphism (of triples).
A similar statement is also true for diagrams with pince-nez, but as
to the former one, one should take also into account the possibility
for each handle to replace its $\lambda$-pince-nez presentation by
the $\mu$-pince-nez one, and vice versa.

A {\it skeleton} of a genus $p$ connected compact oriented
$2$-dimensional manifold with one hole is, by definition, a
semi-skeleton enhanced by a system of pathes that join the marked
point $q\in\partial F$ with the components of the semi-skeleton, the
pathes are called {\it strings}, they are taken disjoint and each of
the $p$ strings is chosen in such a way that in the disc model with
$p$ holes the string riches its  hole at the vertex with outgoing
$\lambda$- and $\mu$-edges.

Now, let us assume that $F$ is equipped with a finite subset
$B\subset F\setminus\partial F$ (later on, such a subset $B$ is
appearing as the branch locus of a finite cover). In such a case, by
a {\it caudate skeleton} we understand a triple $(F,q,S^{\pin})$
where $S^{\pin}$ is a skeleton disjoint from $B$ and extended by a
system of {\it tails}, that is a collection of $n=\vert B\vert$
simple paths connecting the points of $B$ one-by-one with $q$, the
tails being chosen disjoint from each other and from the skeleton.
In particular, $S^\pin$ is homeomorphic to the wedge sum of a
skeleton with $n$ intervals.

\begin{figure}\caption{Plane diagrams of a genus-$2$ surface with its
caudate skeleton.}\label{pinned-picture} \vskip0.1in
\pagestyle{empty}
\newrgbcolor{zzzzzz}{0.6 0.6 0.6}
\newrgbcolor{cccccc}{0.8 0.8 0.8}
\newrgbcolor{xdxdff}{0.49 0.49 1}
\psset{xunit=1.0cm,yunit=1.0cm,algebraic=true,dotstyle=o,dotsize=3pt
0,linewidth=0.8pt,arrowsize=3pt 2,arrowinset=0.25}
\begin{pspicture*}(-3.21,-1.98)(16.45,6.87)
\pscircle(-1.1,3.6){0.99} \psline(-2.09,3.52)(-1.88,3.3)
\psline(-2.09,3.52)(-2.12,3.24) \psline(-1.17,4.59)(-0.84,4.44)
\psline(-0.82,4.72)(-1.17,4.59) \psline(-0.11,3.54)(-0.26,3.36)
\psline(-0.02,3.32)(-0.11,3.54) \psline(-1.44,4.36)(-1.17,4.59)
\psline(-1.52,4.64)(-1.17,4.59) \rput[tl](0.08,3.31){$a_2$}
\rput[tl](0.02,4.58){$b_2$} \rput[tl](-2.57,4.67){$a_2^{-1}$}
\rput[tl](-2.62,3.26){$b_2^{-1}$} \pscircle(6.09,3.55){0.36}
\pscircle(7.43,3.59){0.34} \psline(6.44,3.53)(7.09,3.55)
\psline(6.44,3.53)(6.68,3.39) \psline(6.66,3.71)(6.44,3.53)
\psline(7.77,3.63)(7.55,3.39) \psline(7.77,3.63)(7.55,3.39)
\psline(7.77,3.63)(7.55,3.39) \psline(7.55,3.39)(7.77,3.63)
\psline(7.85,3.32)(7.77,3.63) \rput[tl](7.94,4.06){$a_2$}
\rput[tl](6.62,4.51){$b_2$} \pscircle(0.63,1.18){0.94}
\psline(-0.43,0.91)(-0.31,1.23) \psline(-0.13,0.93)(-0.31,1.23)
\psline(1.39,0.93)(1.57,1.28) \psline(1.67,0.91)(1.57,1.28)
\psline(0.33,1.94)(0.61,2.12) \psline(0.91,1.99)(0.61,2.12)
\psline(0.95,2.26)(0.61,2.12) \psline(0.26,2.22)(0.61,2.12)
\rput[tl](1.36,0.56){$a_1$} \rput[tl](1.47,2.42){$b_1$}
\rput[tl](-0.67,2.42){$a_1^{-1}$} \rput[tl](-0.67,0.7){$b_1^{-1}$}
\pscircle(7.64,1.64){0.37} \pscircle(9,1.62){0.37}
\psline(8.01,1.62)(8.63,1.64) \psline(8.2,1.76)(8.01,1.62)
\psline(8.2,1.5)(8.01,1.62) \psline(9.21,1.44)(9.37,1.65)
\psline(9.48,1.44)(9.37,1.65) \rput[tl](8.16,2.33){$b_1$}
\rput[tl](9.5,2.15){$a_1$} \rput[tl](3.22,3.06){$t_1$}
\rput[tl](2.22,3.47){$t_2$} \rput[tl](10.82,3.24){$t_1$}
\rput[tl](10.03,3.71){$t_2$}
\pscircle[linewidth=1.2pt](0.69,2.92){3.68}
\pscircle[linewidth=1.2pt](8.88,2.9){3.65}
\parametricplot[linewidth=1.2pt]{2.591554445880869}{4.630358903167814}{1*2.22*cos(t)+0*2.22*sin(t)+0.82|0*2.22*cos(t)+1*2.22*sin(t)+1.45}
\parametricplot[linewidth=1.2pt]{-1.2727196395436637}{0.37394794641367346}{1*2.64*cos(t)+0*2.64*sin(t)+-0.13|0*2.64*cos(t)+1*2.64*sin(t)+1.76}
\parametricplot[linewidth=1.2pt]{4.87411180859097}{6.211474740197989}{1*3.21*cos(t)+0*3.21*sin(t)+0.12|0*3.21*cos(t)+1*3.21*sin(t)+2.4}
\parametricplot[linewidth=1.2pt]{-1.273038519549642}{0.11633601375678526}{1*2.84*cos(t)+0*2.84*sin(t)+8.14|0*2.84*cos(t)+1*2.84*sin(t)+1.96}
\parametricplot[linewidth=1.2pt]{-0.8162340983129841}{0.21100514980598906}{1*3.95*cos(t)+0*3.95*sin(t)+6.26|0*3.95*cos(t)+1*3.95*sin(t)+2.12}
\parametricplot[linewidth=1.2pt]{2.6924916576894886}{4.4129820897929894}{1*3.1*cos(t)+0*3.1*sin(t)+9.88|0*3.1*cos(t)+1*3.1*sin(t)+2.21}
\parametricplot[linewidth=1.2pt]{3.0013089405794515}{3.561037801857381}{1*4.38*cos(t)+0*4.38*sin(t)+12.97|0*4.38*cos(t)+1*4.38*sin(t)+1.03}
\psline[linewidth=1.2pt](0.66,0.24)(0.64,-0.76)
\begin{scriptsize}
\psdots[dotstyle=*,linecolor=zzzzzz](-1.07,2.61)
\psdots[dotstyle=*,linecolor=zzzzzz](2.32,2.72)
\psdots[dotstyle=*,linecolor=zzzzzz](-0.11,3.54)
\psdots[dotstyle=*,linecolor=zzzzzz](-1.17,4.59)
\psdots[dotstyle=*,linecolor=zzzzzz](-2.09,3.52)
\psdots[dotstyle=*,linecolor=zzzzzz](0.61,2.12)
\psdots[dotstyle=*,linecolor=zzzzzz](0.66,0.24)
\psdots[dotstyle=*,linecolor=zzzzzz](1.57,1.28)
\psdots[dotstyle=*,linecolor=zzzzzz](-0.31,1.23)
\psdots[dotstyle=*,linecolor=cccccc](7.09,3.55)
\psdots[dotstyle=*,linecolor=zzzzzz](8.63,1.64)
\psdots[dotstyle=*,linecolor=zzzzzz](3.32,2.17)
\psdots[dotstyle=*,linecolor=zzzzzz](10.13,2.95)
\psdots[dotstyle=*,linecolor=zzzzzz](10.96,2.29) \psdots[dotsize=4pt
0,dotstyle=*,linecolor=xdxdff](0.64,-0.76) \psdots[dotsize=4pt
0,dotstyle=*,linecolor=xdxdff](8.97,-0.76)
\end{scriptsize}
\end{pspicture*}
\end{figure}
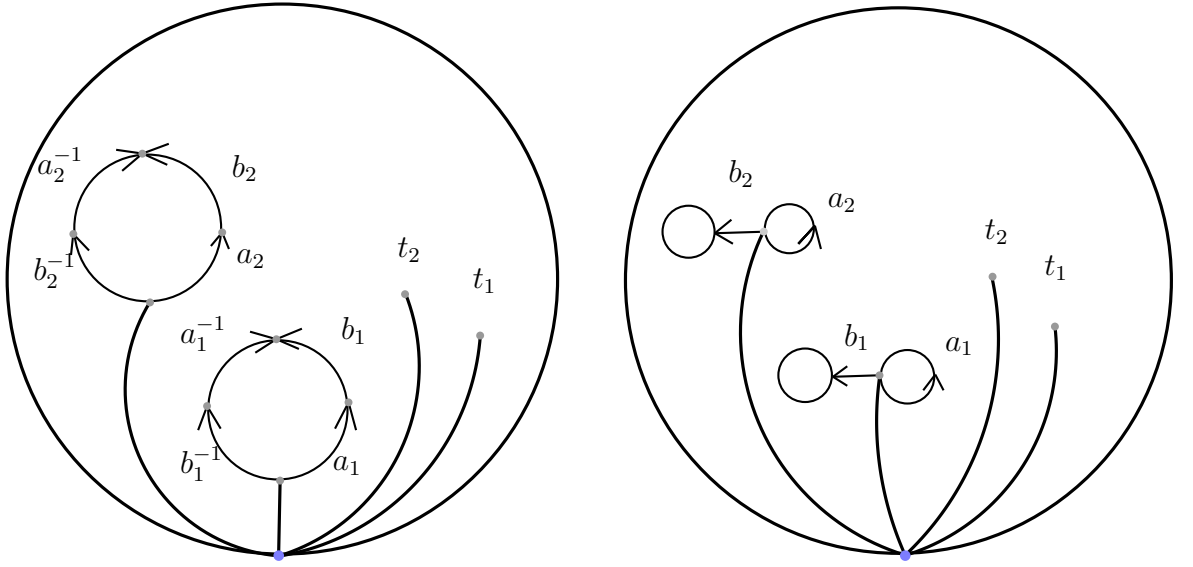

The above notion of open-eyes plane diagrams extends to triples
$(F,q,S)$, where $S$ is either a skeleton or caudate skeleton of
$F$. These diagrams consist of $p=g(F)$ holes or pince-nez in a
disc, a marked point on the boundary of the disc, and a system of
strings, which is enhanced by a system of tails in the case of
caudate skeletons, see Figure \ref{pinned-picture}. Open-eyes plane
diagrams of surfaces with a skeleton or, respectively, caudate
skeleton are defined up to isotopies; and conversely the triple
$(F,q,S)$ is defined up to orientation preserving homeomorphisms of
triples by its open-eyes plane diagram.

In the case of skeletons, and up to isotopy, for a given genus $p$
there is one and only one open-eyes plane diagram. We denote this
diagram by $\Delta_p$ and write $(\Pi_p, q,\Sigma_p)$ to denote the
triple that this diagram defines.

If $p=0$, then $\Delta_0=\Pi_0$ is just the standard disc and
$\Sigma_0$ is reduced to $q$ (the marked point on the boundary of
the disc).

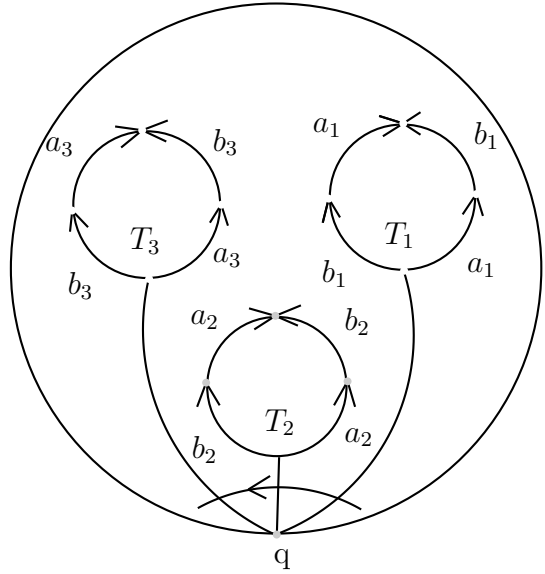
\begin{figure}\caption{
Strings ordered counter-clock wise.}\label{ordering} \vskip0.1in
\pagestyle{empty}
\newrgbcolor{cccccc}{0.8 0.8 0.8}
\psset{xunit=1.0cm,yunit=1.0cm,algebraic=true,dotstyle=o,dotsize=3pt
0,linewidth=0.8pt,arrowsize=3pt 2,arrowinset=0.25}
\begin{pspicture*}(-6.15,-2.66)(12.98,6.77)
\pscircle(0.62,2.75){3.54} \pscircle(-1.1,3.6){0.99}
\pscircle(2.3,3.7){0.98} \psline(-2.09,3.52)(-1.88,3.3)
\psline(-2.09,3.52)(-2.12,3.24) \psline(-1.16,4.59)(-0.84,4.44)
\psline(-0.82,4.72)(-1.16,4.59) \psline(-0.11,3.6)(-0.26,3.36)
\psline(-0.02,3.32)(-0.11,3.6) \psline(-1.44,4.36)(-1.16,4.59)
\psline(-1.52,4.64)(-1.16,4.59) \rput[tl](-0.22,3.02){$a_3$}
\rput[tl](-0.22,4.61){$b_3$} \rput[tl](-2.45,4.49){$a_3$}
\rput[tl](-2.15,2.7){$b_3$} \psline(1.24,3.36)(1.32,3.68)
\psline(1.32,3.68)(1.24,3.36) \psline(1.47,3.43)(1.32,3.68)
\psline(2.04,4.51)(2.31,4.68) \psline(1.99,4.78)(2.31,4.68)
\psline(2.31,4.68)(1.99,4.78) \psline(2.31,4.68)(2.54,4.53)
\psline(2.31,4.68)(2.54,4.83) \psline(3.28,3.74)(3.37,3.45)
\psline(3.28,3.74)(3.1,3.47) \rput[tl](3.16,2.88){$a_1$}
\rput[tl](3.25,4.7){$b_1$} \rput[tl](1.11,4.75){$a_1$}
\rput[tl](1.23,2.86){$b_1$} \psline(6.44,3.53)(6.68,3.39)
\pscircle(0.63,1.18){0.94} \psline(-0.43,0.91)(-0.31,1.23)
\psline(-0.13,0.93)(-0.31,1.23) \psline(1.39,0.93)(1.57,1.25)
\psline(1.67,0.91)(1.57,1.25) \psline(0.33,1.94)(0.61,2.12)
\psline(0.91,1.99)(0.61,2.12) \psline(0.95,2.26)(0.61,2.12)
\psline(0.26,2.22)(0.61,2.12) \rput[tl](1.53,0.63){$a_2$}
\rput[tl](1.53,2.22){$b_2$} \rput[tl](-0.52,2.19){$a_2$}
\rput[tl](-0.5,0.52){$b_2$}
\parametricplot{2.9184591342989785}{4.292997313351126}{1*3*cos(t)+0*3*sin(t)+1.85|0*3*cos(t)+1*3*sin(t)+1.94}
\parametricplot{-1.2059511242595224}{0.3095264935367742}{1*2.83*cos(t)+0*2.83*sin(t)+-0.38|0*2.83*cos(t)+1*2.83*sin(t)+1.86}
\psline(0.63,-0.79)(0.66,0.24)
\parametricplot{1.0464756144914746}{2.078968944630227}{1*2.2*cos(t)+0*2.2*sin(t)+0.65|0*2.2*cos(t)+1*2.2*sin(t)+-2.36}
\psline(0.25,-0.19)(0.54,-0.01) \psline(0.25,-0.19)(0.49,-0.31)
\rput[tl](2.06,3.37){$T_1$} \rput[tl](0.47,0.88){$T_2$}
\rput[tl](-1.32,3.3){$T_3$} \rput[tl](0.59,-0.98){q}
\begin{scriptsize}
\psdots[dotstyle=*,linecolor=white](-1.07,2.61)
\psdots[dotstyle=*,linecolor=white](2.32,2.72)
\psdots[dotstyle=*,linecolor=white](-0.11,3.6)
\psdots[dotstyle=*,linecolor=white](-1.16,4.59)
\psdots[dotstyle=*,linecolor=white](3.28,3.74)
\psdots[dotstyle=*,linecolor=white](1.32,3.68)
\psdots[dotstyle=*,linecolor=white](2.31,4.68)
\psdots[dotstyle=*,linecolor=white](-2.09,3.52)
\psdots[dotstyle=*,linecolor=cccccc](0.61,2.12)
\psdots[dotstyle=*,linecolor=cccccc](1.57,1.25)
\psdots[dotstyle=*,linecolor=cccccc](-0.31,1.23)
\psdots[dotstyle=*,linecolor=cccccc](0.63,-0.79)
\end{scriptsize}
\end{pspicture*}
\end{figure}

If  $p>0$, then the skeleton $\Sigma_p$ contains a non-trivial
semi-skeleton, which we denote by $\Sigma_p^\infty$. A choice of a
skeleton induces in a natural way an ordering on the set of the
components of the semi-skeleton. Namely, we fix the ordering induced
by the counter-clockwise order on the strings of the skeleton, see
Figure \ref{ordering}, and denote the strings of $\Sigma_p$ by
$T_1,\dots, T_p$ following this, counter-clockwise, order. Note that
an ordering of the components of $\Sigma_p^\infty$ is equivalent to
a choice of the skeleton up to isotopy.

Note also that a skeleton $S$ being given, one has a canonical
choice of geometric free generators $\lambda_1,\mu_1, \dots,
\lambda_p,\mu_p$ of the fundamental group $\pi_1(F, q)$, where
$\lambda_i, 1\le i\le p,$ (respectively, $\mu_i$) are represented by
the loops $T_i\star C_{\lambda, i }\star T_i^{-1}$ (respectively,
$T_i\star C_{\mu, i } \star T_i^{-1}$); the numbering respects the
above ordering of the strings.

For given genus $p\ge  1$ and number $n\ge 1$ of tails, the number
of isotopy classes of plane diagrams of surfaces with  caudate
skeletons is greater than $1$ and equal to the binomial coefficient
$C^n_{p+n}$. Indeed, similar to the case of skeletons, a choice of a
caudate skeleton induces a counter-clock wise ordering on the the
set of tails and strings (or, equivalently, on the set that consist
of points of $B$ and the connected components of the semi-skeleton).
Conversely, the counter clock-wise ordering of the set of tails and
strings determines the diagram up to isotopy.

Thus, a caudate skeleton $S^{cut}$ being given, one gets not only a
canonical choice of geometric free generators $\lambda_1,\mu_1,
\dots, \lambda_p,\mu_p$ of the fundamental group $\pi_1(F, q)$, but
also an extension of it to a set of geometric free generators of
$\pi_1(F\setminus B, q)$ by a sequence $\gamma_1,\dots, \gamma_n$
represented by the loops $\Gamma_i\star C_i\star \Gamma_i^{-1}$,
where  $C_i$ denotes a small loop around a point $b_i$ of $B$ and
$\Gamma_i$ a  portion of the tail going to $b_i$. Note that this
whole set of generators of $\pi_1(F\setminus B, q)$ is equipped with
counter clock-wise ordering.

\subsection{Free semigroups of marked coverings}\label{geom-free}
We continuer to consider connected compact oriented $2$-dimensional
manifolds $F$ with one hole and a marked point $q\in \partial F$ and
turn to a study of their ramified finite degree coverings $f:E\to
F$. Let us recall once more that {\it we allow disconnected covering
spaces, but forbid ramifications at the boundary of $F$}. Our aim is
to organize such coverings of a fixed degree $d$ in a semigroup.

To achieve this goal we equip each covering with a marking, that is
a numbering $\nu$ by $1,\dots, d$ of the elements of $f^{-1}(q)$,
and consider the coverings up to certain natural equivalence
relations. Different choices of the equivalence relations lead to
different semigroups.

We start from introducing the {\it free geometric degree $d$
covering semigroup}, which we denote by  $\rm GF\mathbb S_d$. To
build such a semigroup, we equip the base $F$ of each marked
covering $(f : E\to F, \nu, q)$ with a caudate skeleton $S^\pin$
whose tails end at the branch points of the covering. The elements
of $\rm GF\mathbb S_d$ are the triples $(f : E\to F, \nu, S^\pin)$
considered  up to homeomorphisms of coverings respecting all the
ingredients; more precisely, two triples $(f_1 : E_1\to F_1, \nu_1,
S_1^\pin)$ and $(f_2 : E_2\to F_2, \nu_2, S_2^\pin)$ are equivalent
if there are homeomorphisms $\phi :E_1\to E_2$ and $\psi : F_1\to
F_2$ such that $f_2\circ \phi=\psi\circ f_1$,
$\phi\circ\nu_1=\nu_2$, and $\psi(S_1^\pin)=S_2^\pin$.

\begin{figure}
\caption{Plane diagram of a semigroup product.}\label{product}
\vskip0.1in \pagestyle{empty}
\newrgbcolor{xdxdff}{0.49 0.49 1}
\newrgbcolor{zzzzzz}{0.6 0.6 0.6}
\psset{xunit=1.0cm,yunit=1.0cm,algebraic=true,dotstyle=o,dotsize=3pt
0,linewidth=0.8pt,arrowsize=3pt 2,arrowinset=0.25}
\begin{pspicture*}(-4.3,-3.14)(7.52,6.3)
\pscircle[linewidth=1.2pt](-0.52,2.36){2.34}
\pscircle[linewidth=1.2pt](-2.04,2.8){0.42}
\pscircle[linewidth=1.2pt](0.44,2.92){0.42}
\pscircle[linewidth=1.2pt](1.04,1.5){0.36}
\psline[linewidth=1.2pt](-1,2)(-0.52,0.02)
\psline[linewidth=1.2pt](-2.06,2.38)(-0.52,0.02)
\psline[linewidth=1.2pt](-0.52,0.02)(0.36,2.51)
\psline[linewidth=1.2pt](-0.52,0.02)(0.82,1.22)
\psline[linewidth=1.2pt,linestyle=dashed,dash=5pt
5pt](-0.5,4.7)(-0.52,0.02)
\pscircle[linewidth=1.2pt](3.6,2.38){1.13}
\pscircle[linewidth=1.2pt](3.28,2.68){0.4}
\pscircle[linewidth=1.2pt](4,2){0.3}
\psline[linewidth=1.2pt](3.6,1.25)(3.21,2.28)
\psline[linewidth=1.2pt](3.6,1.25)(3.85,1.74)
\pscircle[linewidth=1.2pt](6.22,2.38){1.13}
\pscircle[linewidth=1.2pt](5.9,2.68){0.4}
\psline[linewidth=1.2pt](6.22,1.25)(5.83,2.28)
\psline(6.22,1.25)(6.56,2.18) \rput[tl](2.02,2.6){$=$}
\rput[tl](4.84,2.6){$\bf \cdot$}
\begin{scriptsize}
\psdots[dotsize=4pt 0,dotstyle=*,linecolor=xdxdff](-0.52,0.02)
\psdots[dotsize=4pt 0,dotstyle=*,linecolor=zzzzzz](-1,2)
\psdots[dotsize=4pt 0,dotstyle=*,linecolor=xdxdff](3.6,1.25)
\psdots[dotsize=4pt 0,dotstyle=*,linecolor=xdxdff](6.22,1.25)
\psdots[dotstyle=*,linecolor=zzzzzz](6.56,2.18)
\end{scriptsize}
\end{pspicture*}
\end{figure}
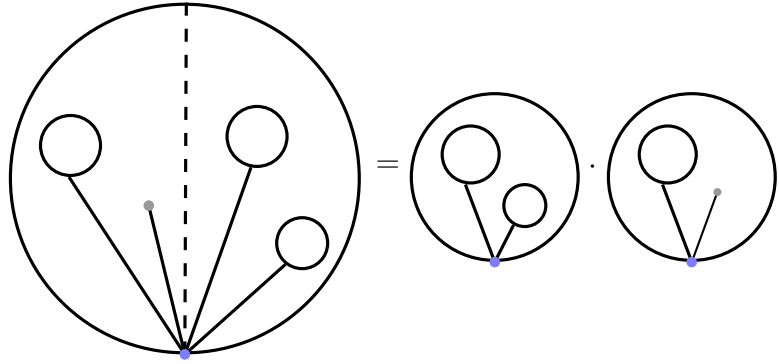

The semigroup structure on  $\rm GF\mathbb S_d$  is defined in a
similar way that was used in \cite{Ku1} in the case of genus zero.
Namely, the product $h=f\cdot g$ of two elements of $\rm GF\mathbb
S_d$ represented by marked ramified coverings $(f:E_1\to F_1, \nu_1,
S_1^\pin)$ and $(g: E_2\to F_2, \nu_2, S_2^\pin)$ (by abuse of
notation we denote by the same symbols both the elements of $\rm
GF\mathbb S_d$ and the underlying coverings) is given by the marked
ramified covering $(h:E\to F, \nu, S_1^{cdt}\cup S_2^{cdt})$, where
$F$ and $h:E\to F$ are obtained, first downstairs, by gluing $F_1$
with $F_2$ along an arc of $\partial F_2$ issued from $q_2$ in the
counter-clockwise direction and an arc of $\partial F_1$ issued from
$q_1$ in the clockwise direction, (see Figure \ref{product}) and,
second upstairs, by a gluing of $f$ and $g$ that preserves  the
markings over $q=q_1=q_2$. This operation respects the equivalence
relation.

We equip $\rm GF\mathbb S_d$ with a map $\alpha : {\rm GF\mathbb
S_d} \to \mathcal S_d$ that evaluates the boundary monodromy (taken
in the direction of boundary orientation). As it follows from the
gluing procedure, this map is a homomorphism. Furthermore, the
symmetric group $\mathcal S_d$ naturally acts on $\rm GF\mathbb S_d$
by renumbering the points of the fibre $f^{-1}(q)$. Thus, $\rm
GF\mathbb S_d$ becomes in a canonical way a semigroup over $\mathcal
S_d$.

For each $g\in \mathcal S_d$ denote by $X_g\in \rm GF\mathbb S_d$
the element represented by a ramified covering  $f: E\to \Pi_0$ with
one branch point, a marked point $q\in\partial\Pi_0$, and the
monodromy $\alpha(X_g)$ equal to $g$. Such an element is defined
uniquely, as it follows, for example, from Proposition
\ref{completion}.

Next, consider the torus $\Pi_1$ with a hole, a marked point
$q\in\partial\Pi_1$ and a skeleton $\Sigma_1$ that includes the
semi-skeleton $\Sigma^\infty_1$ and the string $T_1$. Pick a pair
$a,b\in \mathcal S_d$ and denote by $Y_{a,b}\in \rm GF\mathbb S_d$
the element represented by a non-ramified covering $f:E\to \Pi_1$
with monodromy $a$ along the loop $T_1\star C_\lambda \star
T_1^{-1}$ and $b$ along the loop $T_1\star C_\mu \star T_1^{-1}$.
Such an element is also defined uniquely, as it follows, for
example, from  Proposition \ref{encoding}. Note that
$\alpha(Y_{a,b})=aba^{-1}b^{-1}$.

\begin{prop}\label{geom-norel}  The semigroup $\rm GF\mathbb S_d$ is a free semigroup over the group $\mathcal S_d$, its set of free
generators is formed by $X_g, g\in \mathcal S_d\setminus \{\bf 1\}$,
and $Y_{a,b}, a,b\in \mathcal S_d$.
\end{prop}
\proof Follows from the following isotopy unicity:  for a  surface
with a given caudate skeleton, its open-eyes caudate skeleton plane
diagram $(\Sigma_p^\pin,q)\subset (\Delta_p,q)$  is unique up to
isotopies in $(\Delta_p,q)$. \qed

As a set, the semigroup $\rm GF\mathbb S_d$  splits  in a disjoint
union of subsets, $(\rm GF\mathbb S_d)_{n,p}$, that correspond to
coverings with a given number $n$ of branch points over surfaces of
given genus $p$, or, saying in another way, to words with $n$
letters $X_g$ and $p$ letters $Y_{a,b}$.

\subsection{Very strong semigroup of marked coverings}\label{geom-frozen} As a
next step, we replace caudate skeletons by skeletons, that is forget
the tails going to the branch points, and thus construct another
semigroup over $\mathcal S_d$ replacing everywhere in the above
construction of $\rm GF\mathbb S_d$ the caudate skeletons by
skeletons. We call this new semigroup the {\it very strong semigroup
of degree $d$ marked coverings} and denote it by $\rm GV\mathbb
S_d$. The forgetful map $\rm GF\mathbb S_d\to \rm GV\mathbb S_d$
consisting in replacing a caudate skeleton by the skeleton is a well
defined homomorphism of semigroups over $\mathcal S_d$. Let us\marginpar{Kh: us}
denote
by the same symbols $X_g, Y_{a,b}$ the images in $\rm GV\mathbb S_d$
of the above free generators $X_g, Y_{a,b}\in \rm GF\mathbb S_d$.

\begin{prop}\label{frozen-rel}  The elements $X_g, g\in \mathcal S_d\setminus\{\bf 1\},$ and
$Y_{a,b}, a,b\in \mathcal S_d,$ form a set of generators of the
semigroup  $\rm GV\mathbb S_d$. They satisfy the relations
\begin{equation} \label{fr12'} X_{g}\cdot
Y_{a,b}=Y_{a,b}\cdot X_{g^{[a,b]}}
\end{equation} and
\begin{equation} \label{fr11} X_{g_1}\cdot
X_{g_2}=X_{g_2}\cdot\,X_{g_1^{g_2}}
\end{equation}
for any ${g_1}, {g_2},a,b\in \mathcal S_d$. These are the defining
relations of  $\rm GV\mathbb S_d$.
\end{prop}

\begin{figure}\caption{Plane diagram of an elementary change of ordering.}\label{Fig-ordering}
\newrgbcolor{zzzzzz}{0.6 0.6 0.6}
\newrgbcolor{xdxdff}{0.49 0.49 1}
\psset{xunit=1.0cm,yunit=1.0cm,algebraic=true,dotstyle=o,dotsize=3pt
0,linewidth=0.8pt,arrowsize=3pt 2,arrowinset=0.25}
\begin{pspicture*}(-4.3,-0.42)(8.5,6.3)
\pscircle[linewidth=1.2pt](-0.82,3.2){2.2} \rput[tl](2.14,3.7){$ =
$} \pscircle[linewidth=1.2pt](5.24,3.2){2.16}
\pscircle[linewidth=1.2pt](-0.84,1.98){0.51}
\pscircle[linewidth=1.2pt](5.26,1.98){0.5}
\parametricplot[linewidth=1.2pt,linestyle=dotted]{-1.3391635018811003}{1.2131915747686608}{1*1.15*cos(t)+0*1.15*sin(t)+-1.08|0*1.15*cos(t)+1*1.15*sin(t)+4.32}
\parametricplot[linewidth=1.2pt,linestyle=dotted]{1.657665644272418}{4.735332946487752}{1*1.1*cos(t)+0*1.1*sin(t)+-0.72|0*1.1*cos(t)+1*1.1*sin(t)+2.11}
\parametricplot[linewidth=1.2pt,linestyle=dotted]{1.828718819819892}{4.344959636305275}{1*1.14*cos(t)+0*1.14*sin(t)+5.65|0*1.14*cos(t)+1*1.14*sin(t)+4.26}
\parametricplot[linewidth=1.2pt,linestyle=dotted]{-1.2410504917360088}{1.3330484787095214}{1*1.13*cos(t)+0*1.13*sin(t)+4.97|0*1.13*cos(t)+1*1.13*sin(t)+2.1}
\parametricplot[linewidth=1.2pt,linestyle=dashed,dash=2pt 2pt]{-1.4578366064746202}{1.4698844161994373}{1*1.63*cos(t)+0*1.63*sin(t)+5.16|0*1.63*cos(t)+1*1.63*sin(t)+2.66}
\parametricplot[linewidth=1.2pt,linestyle=dashed,dash=2pt 2pt]{1.5396803794493856}{4.767821641829171}{1*1.67*cos(t)+0*1.67*sin(t)+-0.79|0*1.67*cos(t)+1*1.67*sin(t)+2.67}
\psline[linewidth=1.2pt](-0.7,1.01)(-0.71,1.49)
\psline[linewidth=1.2pt](5.34,1.04)(5.35,1.49)
\begin{scriptsize}
\psdots[dotsize=4pt 0,dotstyle=*,linecolor=zzzzzz](-0.74,4.34)
\psdots[dotsize=4pt 0,dotstyle=*,linecolor=zzzzzz](5.32,4.28)
\psdots[dotsize=4pt 0,dotstyle=*,linecolor=xdxdff](-0.7,1.01)
\psdots[dotsize=4pt 0,dotstyle=*,linecolor=xdxdff](5.34,1.04)
\end{scriptsize}
\end{pspicture*}
\end{figure}
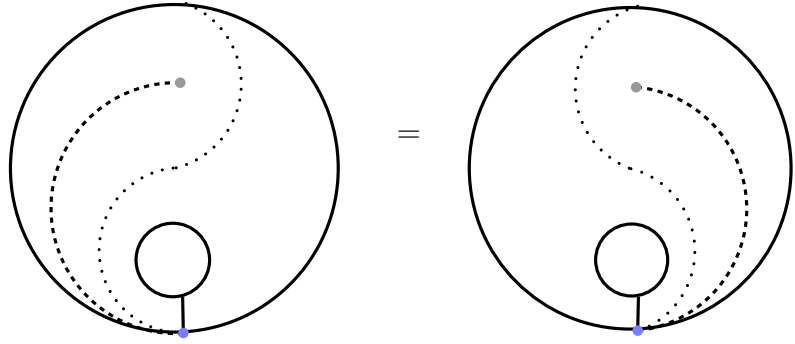

\proof For a surface with a given skeleton and branch locus, an
enhancing of its open-eyes skeleton plane diagram
$(\Sigma_p,q)\subset( \Delta_p,q)$ to an open-eyes caudate skeleton
diagram $(\Sigma^\pin_p,q)\subset( \Delta_p,q)$ with a given
ordering of strings and tails is unique up to isotopies of
$(\Sigma_p\cup B,q)$ in $(\Delta_p,q)$. The relation (\ref{fr12'})
reflects the elementary change of ordering (see Figure
\ref{Fig-ordering}), and (\ref{fr11}) the standard Artin-Hurwitz
half-twist of two points of $B$ in $\Delta_p\setminus\Sigma_p$ (cf.,
\cite{Ku3}). \qed

Similar to  $\rm GF\mathbb S_d$, the semigroup $\rm GV\mathbb S_d$
splits as a set in a disjoint union of subsets, $(\rm GV\mathbb
S_d)_{n,p}$, that correspond to coverings with a given number $n$ of
branch points over surfaces of given genus $p$.

\subsection{Strong versal geometric covering semigroups}\label{geom-strong}
These semigroups are most close to the classical theory of Hurwitz
spaces.

We equip the base of each covering with a skeleton and consider the
triples $(f_1 : E_1\to F_1, \nu_1, S_1)$ and $(f_2 : E_2\to F_2,
\nu_2, S_2)$, where $S_1$ and $S_2$ are skeletons of $F_1$ and
$F_2$, respectively, up to the equivalence generated by two binary
relations: first, up to isotopy of coverings with fixed base,
marking, and skeleton, but moving branch points (which can move, in
particular, through the skeleton); second, up to homeomorphisms
respecting all the ingredients, that is up to homeomorphisms $\phi
:E_1\to E_2$ and $\psi : F_1\to F_2$ such that $f_2\circ
\phi=\psi\circ f_1$, $\phi\circ\nu_1=\nu_2$, and
$\psi(S_1^\infty)=S_2^\infty$. The only, but major, difference with
respect to the previous very strong covering semigroups is that we
authorize the branch points to cross the skeleton.

Taking into account this additional equivalence relation, we obtain
another semigroup over $\mathcal S_d$, which we denote  by $\rm
G\mathbb S_d$ and call the  {\it strong versal geometric degree $d$
covering semigroup}. The quotient map $\rm GV\mathbb S_d\to \rm
G\mathbb S_d$ is a homomorphism of semigroups over $\mathcal S_d$
and, set theoretically, it splits in quotient maps $(\rm GV\mathbb
S_d)_{n,p}\to( \rm G\mathbb S_d)_{n,p}$.

For each $n$ and $p$, let us fix the surface $F$ (of genus $p$) and
its skeleton $S$, and place the branch locus $B$ (of cardinality
$n$) to be disjoint from $S$. Then the braid group on $n$ strands,
that is the group $Br_n(F,\partial F)$ of isotopy classes of
orientation preserving identical on the boundary self-homeomorphisms
of $(F,B)$, becomes to act naturally on $(\rm GV\mathbb S_d)_{n,p}$,
and, as it follows also directly from the definitions, the fibers of
the quotient map $(\rm GV\mathbb S_d)_{n,p}\to( \rm G\mathbb
S_d)_{n,p}$ are the orbits of this action.

Let us recall that, on the other hand, the braid group
$Br_n(F,\partial F)$ can be canonically identified with the
fundamental group $\pi_1(F^{(n)}\setminus \Delta)$, where $F^{(n)}$
is the symmetric product of $n$ copies of $F$ and $\Delta$ is the
discriminant locus, that is, the set of those $n$-tuples that
contain fewer than $n$ distinct points. More precisely, we start
from fixing a set $B=\{P_1,\dots,P_n\}\subset F\setminus
\partial F$ consisting of $n$ distinct points and treat $Br_n(F,\partial F)$ as
the group of homotopy classes of geometric braids, where as is
usual: by a {\it geometric braid} on $F$ based at $B$ we understand
an $n$-tuple $\Psi= (\psi_1,\dots, \psi_n)$ of paths
$\psi_i:[0,1]\to F\setminus\partial F$ such that
\begin{itemize}
\item[$(1)$] $\psi_i(0)=P_i$ and $\psi_i(1)\in B$ for each
$i=1,\dots,n$;
\item[$(2)$] $\psi_1(t),\dots,\psi_n(t)$ are distinct points of $F\setminus\partial F$
for each $t\in [0,1]$;
\end{itemize}
and multiplication is given by concatenation of paths.

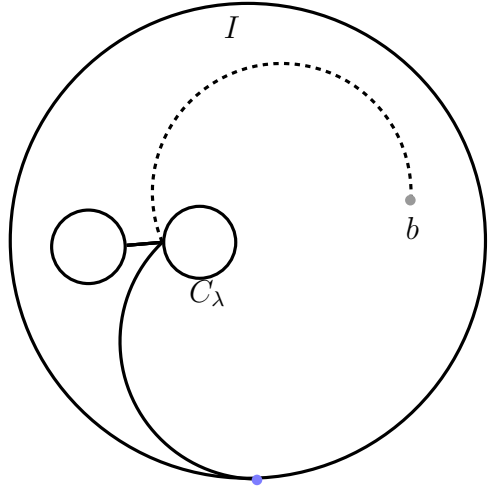
\begin{figure}
\caption{Simple path joining a branch point with a bouquet
vertex.}\label{simple-path}
\newrgbcolor{xdxdff}{0.49 0.49 1}
\newrgbcolor{zzzzzz}{0.6 0.6 0.6}
\psset{xunit=1.0cm,yunit=1.0cm,algebraic=true,dotstyle=o,dotsize=3pt
0,linewidth=0.8pt,arrowsize=3pt 2,arrowinset=0.25}
\begin{pspicture*}(-5.98,-1.38)(11.34,6.42)
\pscircle[linewidth=1.2pt](0.68,2.58){3.18}
\pscircle[linewidth=1.2pt](-1.44,2.5){0.51}
\pscircle[linewidth=1.2pt](0.04,2.56){0.5}
\psline[linewidth=1.2pt](-0.93,2.52)(-0.46,2.56)
\psline[linewidth=1.2pt](-0.93,2.52)(-0.46,2.56)
\psline[linewidth=1.2pt](-0.46,2.56)(-0.93,2.52)
\parametricplot[linewidth=1.2pt]{2.3360229059375843}{4.706656358741243}{1*1.83*cos(t)+0*1.83*sin(t)+0.81|0*1.83*cos(t)+1*1.83*sin(t)+1.24}
\parametricplot[linewidth=1.2pt,linestyle=dashed,dash=2pt 2pt]{-0.05686189859425994}{3.5346058331201746}{1*1.72*cos(t)+0*1.72*sin(t)+1.13|0*1.72*cos(t)+1*1.72*sin(t)+3.22}
\rput[tl](2.78,2.9){$ b $} \rput[tl](-0.1,2.06){$ C_\lambda $}
\rput[tl](0.36,5.56){$  I $}
\begin{scriptsize}
\psdots[dotsize=4pt 0,dotstyle=*,linecolor=xdxdff](0.8,-0.6)
\psdots[dotsize=4pt 0,dotstyle=*,linecolor=zzzzzz](2.84,3.12)
\end{scriptsize}
\end{pspicture*}
\end{figure}

By  a $\lambda$- (respectively, $\mu$-) move we understand a
geometric braid whose all but one strands are constant and the
remaining one follows a path $ I \star C_\lambda \star I^{-1}$
(respectively, $ I \star C_\mu \star I^{-1}$) where $I$ is a simple
path in the complement of $S\cup (B\setminus \{ b\})$ joining a
point $b\in B$ with the vertex of a bouquet $C_\lambda\vee
C_\mu\subset S$, see Figure \ref{simple-path}. The standard
Artin-Hurwitz (half-twist) geometric braids exchanging two points of
$B$ will be called $H$-moves.

The following proposition is well known. In fact, it follows easily,
for example, from the exact sequences $$ 1\to \pi_1(F\setminus
B')\to PBr_n(F,\partial F)\to PBr_{n-1}(F,\partial F)\to 1
$$ and $$ 1\to PBr_n(F,\partial F)\to Br_n(F,\partial F)\to \mathcal
S_n\to 1, $$ where the second sequence is the definition of the pure
braid groups, $PBr_n$, and $B'$ denotes $B\setminus \{P_1\}$.

\begin{prop}\label{generators-braid} {\rm (\cite{FaN,FoN})} The braid group $Br_n(F,\partial F)$
is generated by $H$-, $\lambda$-, and $\mu$-moves.\qed
\end{prop}

The following claim is an immediate consequence of the above
Proposition.

\begin{cor}\label{2-terms} The relations imposed by the partial quotient maps $(\rm GV\mathbb S_d)_{2,0}\to (\rm G\mathbb S_d)_{2,0}$
and $(\rm GV\mathbb S_d)_{1,1}\to (\rm G\mathbb S_d)_{1,1}$ imply
all the other  relations imposed by the quotient map $\rm GV\mathbb
S_d\to \rm G\mathbb S_d$.\qed
\end{cor}

To describe finally the semigroup $G\mathbb S_d$ in terms of
generators and relations, let us denote by the same symbols $X_g,
Y_{a,b}$ the images in $ \rm G\mathbb S_d$ of $X_g, Y_{a,b}\in \rm
GV\mathbb S_d $.

\begin{prop}\label{geom-rel} The elements $X_g, g\in \mathcal S_d\setminus\{\bf 1\},$ and
$Y_{a,b}, a,b\in \mathcal S_d,$ form a set of generators of the
semigroup $\rm G\mathbb S_d$. They satisfy the relations:
\begin{equation} \label{grel11} X_{g_1}\cdot
X_{g_2}=X_{g_2}\cdot\,X_{g_1^{g_2}}
\end{equation}
for any ${g_1}, {g_2}\in \mathcal S_d$; and
\begin{equation} \label{grel12'} X_{g}\cdot
Y_{a,b}=Y_{a,b}\cdot X_{g^{[a,b]}},
\end{equation}
\begin{equation} \label{grel13'} X_{g}\cdot
Y_{a,b}=X_{g^{c_1}}\cdot Y_{ga,b}, \quad c_1=ab^{-1}a^{-1}g^{-1},
\end{equation}
\begin{equation} \label{grel14'}
Y_{a,b} \cdot X_{g}=Y_{a,g^{-1}b}\cdot X_{g^{c_3}}, \quad
c_2=ba^{-1}b^{-1}g,
\end{equation}
for any ${g}, a,b\in \mathcal S_d$. These are the defining relations
of $\rm G\mathbb S_d$.
\end{prop}

\proof Due to Propositions  \ref{frozen-rel}, \ref{generators-braid}
and Corollary \ref{2-terms}, the only new, with respect to
Proposition \ref{frozen-rel}, relations are given by the  $\lambda$-
and  $\mu$-moves  in $(\rm GV\mathbb S_d)_{1,1}$.

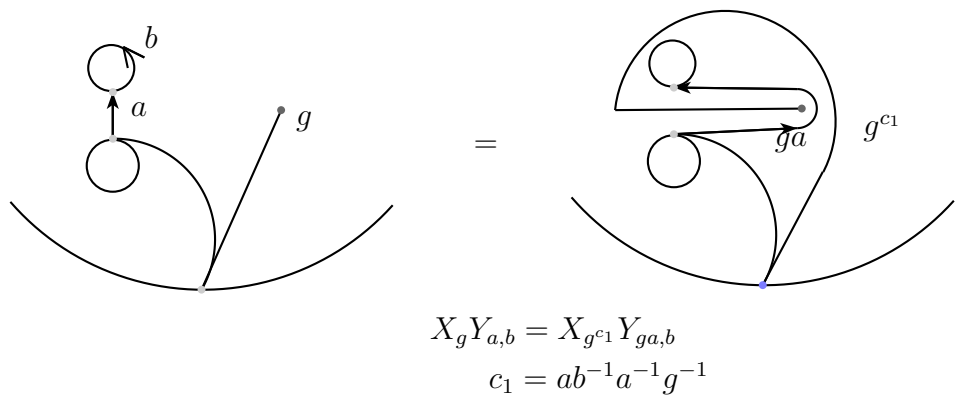
\begin{figure}\caption{Moving a branch point through an $\lambda$-cut in an $XY$-product}\label{Fa}
\vskip0.1in \pagestyle{empty}
\newrgbcolor{cccccc}{0.8 0.8 0.8}
\newrgbcolor{wwwwww}{0.4 0.4 0.4}
\newrgbcolor{xdxdff}{0.49 0.49 1}
\psset{xunit=1.0cm,yunit=1.0cm,algebraic=true,dotstyle=o,dotsize=3pt
0,linewidth=0.8pt,arrowsize=3pt 2,arrowinset=0.25}
\begin{pspicture*}(-4.3,-1.86)(13.8,6.3)
\parametricplot{3.8559807057465623}{5.55123413606367}{1*3.39*cos(t)+0*3.39*sin(t)+-1.28|0*3.39*cos(t)+1*3.39*sin(t)+5.96}
\pscircle(-2.48,4.22){0.36} \pscircle(-2.5,5.52){0.32}
\psline(-2.48,5.2)(-2.48,4.58)
\parametricplot{-0.5263728493588333}{1.5887171483999096}{1*1.34*cos(t)+0*1.34*sin(t)+-2.46|0*1.34*cos(t)+1*1.34*sin(t)+3.24}
\psline(-0.24,4.96)(-1.3,2.57)
\parametricplot{3.8559807057465623}{5.55123413606367}{1*3.39*cos(t)+0*3.39*sin(t)+6.18|0*3.39*cos(t)+1*3.39*sin(t)+6.02}
\pscircle(4.98,4.28){0.36} \pscircle(4.96,5.58){0.32}
\parametricplot{-0.5263728493588324}{1.588717148399909}{1*1.34*cos(t)+0*1.34*sin(t)+5|0*1.34*cos(t)+1*1.34*sin(t)+3.3}
\parametricplot{-1.5707963267948966}{1.5707963267948966}{1*0.26*cos(t)+0*0.26*sin(t)+6.62|0*0.26*cos(t)+1*0.26*sin(t)+4.98}
\psline(4.98,4.64)(6.62,4.72) \psline(4.98,5.26)(6.62,5.24)
\psline(6.68,4.98)(4.2,4.96)
\parametricplot{-0.4755203853156553}{3.0395221253407603}{1*1.47*cos(t)+0*1.47*sin(t)+5.66|0*1.47*cos(t)+1*1.47*sin(t)+4.81}
\psline(6.16,2.63)(6.96,4.14) \psline{->}(-2.48,4.58)(-2.48,5.2)
\psline{->}(4.98,4.64)(6.62,4.72) \psline{->}(6.62,5.24)(4.98,5.26)
\psline(-2.34,5.8)(-2.06,5.66) \psline(-2.34,5.8)(-2.06,5.66)
\psline(-2.28,5.52)(-2.34,5.8) \rput[tl](-2.24,5.06){$a$}
\rput[tl](-2.06,6.08){$b$} \rput[tl](-0.04,4.94){$g$}
\rput[tl](2.32,4.56){=} \rput[tl](6.32,4.66){$ga$}
\rput[tl](7.5,4.9){$g^{c_1}$} \rput[lt](1.74,2.2)
{$X_gY_{a,b}=X_{g^{c_1}}Y_{ga,b}$}
\rput[lt](2.38,1.6){ $c_1=ab^{-1}a^{-1}g^{-1}$}
\begin{scriptsize}
\psdots[dotstyle=*,linecolor=cccccc](-1.3,2.57)
\psdots[dotstyle=*,linecolor=cccccc](-2.48,4.58)
\psdots[dotstyle=*,linecolor=cccccc](-2.48,5.2)
\psdots[dotstyle=*,linecolor=wwwwww](-0.24,4.96)
\psdots[dotstyle=*,linecolor=xdxdff](6.16,2.63)
\psdots[dotstyle=*,linecolor=cccccc](4.98,4.64)
\psdots[dotstyle=*,linecolor=cccccc](4.98,5.26)
\psdots[dotstyle=*,linecolor=wwwwww](6.68,4.98)
\end{scriptsize}
\end{pspicture*}
\end{figure}

\begin{figure}\caption{Moving a branch point through a $\mu$-cut in a $YX$-product}\label{Fc}
\vskip0.1in \pagestyle{empty}
\newrgbcolor{xdxdff}{0.49 0.49 1}
\newrgbcolor{cccccc}{0.8 0.8 0.8}
\psset{xunit=1.0cm,yunit=1.0cm,algebraic=true,dotstyle=o,dotsize=3pt
0,linewidth=0.8pt,arrowsize=3pt 2,arrowinset=0.25}
\begin{pspicture*}(-4.3,-2.52)(19.82,6.3)
\parametricplot{4.124386376837123}{5.2922025040268785}{1*5.55*cos(t)+0*5.55*sin(t)+6.72|0*5.55*cos(t)+1*5.55*sin(t)+5.24}
\pscircle(6.4,2.4){0.38} \pscircle(4.78,2.38){0.38}
\psline(6.02,2.36)(6.99,-0.31) \rput[tl](6.96,2.62){$a$}
\rput[tl](4.52,0.9){$b'=g^{-1}b$} \rput[tl](5.44,1.96){$g^{c_2}$}
\parametricplot{4.140437878210407}{5.302791579324273}{1*5.88*cos(t)+0*5.88*sin(t)+-0.66|0*5.88*cos(t)+1*5.88*sin(t)+5.52}
\psline(-0.54,-0.35)(-1.3,1.4) \pscircle(-0.62,2.26){0.36}
\psline(-0.98,2.26)(-1.48,2.26) \pscircle(-1.84,2.26){0.36}
\psline{->}(-0.98,2.26)(-1.48,2.26) \psline(-0.98,2.26)(-0.54,-0.35)
\rput[tl](-1.44,0.66){g} \rput[tl](0,2.5){a}
\rput[tl](-1.22,2.94){b} \psline(-0.34,2.3)(-0.34,2.48)
\psline(-0.18,2.38)(-0.34,2.48) \rput[tl](3.1,1.96){=}
\rput[lt](3.06,-0.72)
{$Y_{a,b}X_g=Y_{a,g^{-1}b}X_{g^{c_2}}$}
\rput[lt](3.88,-1.6) {$c_2=ba^{-1}b^{-1}g$}
\parametricplot{3.113822016996372}{6.255414670586165}{1*0.36*cos(t)+0*0.36*sin(t)+5.64|0*0.36*cos(t)+1*0.36*sin(t)+1.39}
\psline(6.02,2.36)(6,1.38) \psline{->}(5.28,1.4)(5.16,2.42)
\begin{scriptsize}
\psdots[dotstyle=*,linecolor=xdxdff](6.99,-0.31)
\psdots[dotstyle=*,linecolor=cccccc](5.68,1.44)
\psdots[dotstyle=*,linecolor=xdxdff](-0.54,-0.35)
\psdots[dotstyle=*,linecolor=cccccc](-1.3,1.4)
\end{scriptsize}
\end{pspicture*}
\end{figure}
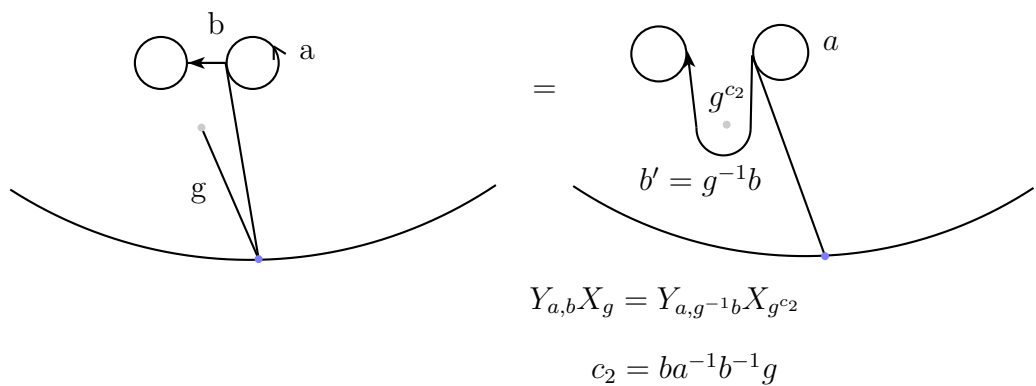

Under pulling the branch point through  $\lambda$- or $\mu$-cuts as
it is shown in Figures \ref{Fa} and \ref{Fc}, we obtain the the
relations (\ref{grel13'}) and (\ref{grel14'}). In particular,  the
elements $g'$ appearing in the relations $X_g\cdot
Y_{a,b}=X_{g'}\cdot Y_{ga,b}$ and, respectively, $Y_{a,b}\cdot
X_g=Y_{a,g^{-1}b}\cdot X_{g'}$ (see Figures \ref{Fa} and \ref{Fc})
can be found, for example, from the identities expressing the
unchanged monodromy along the hole:
$gaba^{-1}b^{-1}=g'a'b(a')^{-1}b^{-1}, a'=ga$ under the move across
the $\lambda$-cut, as is shown in Figure \ref{Fa}, and
$gaba^{-1}b^{-1}=g'ab'a^{-1}(b')^{-1}, b'=g^{-1}b$ under the
$\mu$-cut, as is shown in Figure \ref{Fc}. \qed

\begin{cor}\label{add-rel}
The elements $X_g, g\in \mathcal S_d\setminus\{\bf 1\},$ and
$Y_{a,b}, a,b\in \mathcal S_d,$ satisfy the relation
\begin{equation} \label{grel19}
X_g\cdot Y_{a,b}=X_{g^{[a,b]}}\cdot Y_{a^{g^{[a,b]}},b^{g^{[a,b]}}}
\end{equation}
\end{cor}

\proof This relation follows from the relations (\ref{grel12'}),
(\ref{grel13'}), and (\ref{grel14'}). It expresses the braiding of a
branch point around the handle (that is around the hole on the
open-eyes plane diagram) . \qed

Due to Proposition \ref{geom-rel} and Corollary \ref{add-rel}, there
arises a canonical morphism
$$
\varrho_{d } : \rm G\mathbb S_d\to \rm \mathbb S_d
$$
that maps $X_g$ to $x_g$ and $Y_{a,b}$ to $y_{a,b}$, see subsection
\ref{def-cov}.

\begin{prop}\label{geom=alg} The morphism $\varrho_{d
}: \rm G\mathbb S_d\to \rm \mathbb S_d$ is an isomorphism of
semigroups over $\mathcal S_d$.
\end{prop}
\proof The morphism $\varrho_{d}$ literally translates the list of
generators and defining relations in $ \rm G\mathbb S_d$ given in
Proposition \ref{geom-rel}  and extended by ({\ref{grel19}) into the
list of generators and defining relations in $\mathbb S_d$ given in
Proposition \ref{red-rel}. \qed

\subsection{Admissible geometric covering
semigroups}\label{weak-semigroups} Again, to  define new geometric
covering semigroups, we start  from the very strong semigroup of
degree $d$ marked coverings, $\rm GV\mathbb S_d$, and add
supplementary relations between the triples $(f : E\to F, \nu,
S^{cdt})$ representing the elements of $\rm GV\mathbb S_d$. For that
purpose, it is convenient to introduce some auxiliary category
$\mathcal F$.

The objects of $\mathcal F$ are the triples $(F,q,S^{cdt})$ where
$F$ is a connected compact oriented surface with one hole, $q$ is a
fixed point on its boundary, and $S^{cut}$ is a caudate skeleton of
$F$; we denote the set of ends of the tails of $S^{cdt}$ by $B$, its
cardinality by $n$, and the genus of $F$ by $p$. The morphisms of
$\mathcal F$ are the preserving orientation homeomorphisms of
triples $(F,q,B)$.\footnote{Note that morphisms are not supposed to
respect skeletons.}

We denote by $Homeo$ the whole set of morphisms of $\mathcal F$; by
$Homeo_{n,p}$ the subset of  $Homeo$ consisting of the morphisms
between the  triples $(F,q,S^{cdt})$ with given  $|B|=n$ and
$g(F)=p$; and  by $H(F,q,S^{cdt})$ the group consisting of the
self-homeomorphisms $(F,q,B)\to (F,q,B)$.\footnote{Note that
$H(F,q,S^{cut})$ is in a canonical bijection with the set of
self-morphisms $(F,q,S^{cdt})\to (F,q,S^{cdt})$.} We say that a
collection of subsets $\widetilde H_{n,p}\subset Homeo_{n,p}$,
$n\geq 0$, $p\geq 0$, is {\it geometrically admissible}, if it
contains the isotopies of $B$ in $F\setminus \partial F$ and for
each two triples $(F,q,S^{cdt})$ and $(F,q,S^{cdt})$ with the same
$n$ and $p$ there is a morphism $(F,q,S^{cdt})\to (F,q,S^{cdt})$
belonging to $\widetilde H_{n,p}$.

For each $n$ and $p$, we fix  the triple $(F,q,S^{cdt})$. The
caudate skeleton $S^{cdt}$ defines a frame in the free group
$\pi_1(F\setminus B,q)$, while each element $\psi\in H(F,q,S^{cdt})$
defines an automorphism of this free group  $\mathbb
F^{n+2p}=\pi_1(F\setminus B,q)$. Let us denote by $H_{n,p}$ the
subgroups of  $Aut(\mathbb F^{n+2p})$, $\mathbb
F^{n+2p}=\pi_1(F\setminus B,q)$, representing $\widetilde
H_{n,p}\cap H(F,q,S^{cut}) $.

As follows directly from the definitions,  $H_{n,p}$
 are admissible  subgroups of  $Aut(\mathbb F^{n+2p})$, if the
 collection $\widetilde H_{n,p}\subset Homeo_{n,p}$ is geometrically admissible.
For example, one get geometrically admissible collections by
considering the homeomorphisms preserving $\lambda$-circles up to
isotopy (respectively, $\mu$-circles), or taking the whole sets,
$\widetilde H_{n,p}=Homeo_{n,p}$. In fact, one can easily go other
way round and, starting from a collection of admissible subgroups
$H_{n,p}$ of $Aut(\mathbb F^{n+2p})$, build a geometrically
admissible collection by attributing to $\tilde H_{n,p}\cap
H(F,q,S^{cdt})$ all the elements in $H(F,q,S^{cdt})$ that act in
$\mathbb F^{n+2p}=\pi_1(F\setminus B,q)$ as elements of $H_{n,p}$
(with respect to the frame defined by $S^{cdt}$).

Let $\mathcal {\tilde H}=\{ \widetilde H_{n,p}\}_{n\geq 0,p\geq0}$
be a geometrically admissible collection. Assume in addition that
this collection of homeomorphisms is closed under the boundary
connected sum of the triples $(F,q,B)$ (see definition of the
product in the semigroup $\rm GF\mathbb S_d$). We say that two
triples $(f_1 : E_1\to F_1, \nu_1, S^{cdt}_1)$ and $(f_2 : E_2\to
F_2, \nu_2, S^{cdt}_2)$ are $\mathcal {\tilde H}$-{\it equivalent}
if there are homeomorphisms $\phi :E_1\to E_2$ and $\psi :  (F_1,
q_1, B_1) \to (F_2,q_2, B_2)$
 such that $f_2\circ \phi=\psi\circ f_1$,
$\phi\circ\nu_1=\nu_2$ , and $\psi \in \widetilde H_{n,p}$. By means
of such an additional equivalence relation, we obtain a semigroup
over the group $\mathcal S_d$ taking the quotient $\rm G\mathbb
S_{d,\mathcal {\tilde H}-equiv}=GV\mathbb S_d/\{
f_1\stackrel{\mathcal {\tilde H}}{\sim}f_2\}$; we call the
semigroups thus obtained {\it admissible versal geometric degree $d$
covering semigroups}. In particular, if $\widetilde
H_{n,p}=Homeo_{n,p}$ for each $n$ and $p$, the semigroup $\rm
G\mathbb S_{d,\mathcal {\tilde H}-equiv}$ is called the {\it weak
versal geometric degree $d$ covering semigroup} and denoted by
$GW\mathbb S_{d}$.

Let $G$ be a subgroup of $\mathcal S_d$ and $O\subset G$ be some its
equipment.  A subsemigroup ${\rm G\mathbb S_{d,\mathcal {\tilde
H}-equiv}}(G,O)$ of the semigroup ${\rm G\mathbb S_{d,\mathcal
{\tilde H}-equiv}}$ generated by $X_g$, $g\in O$, and $Y_{a,b}$,
$a,b\in G$, is a semigroup over $G$. Its elements are $\mathcal
H$-equivalence classes of degree $d$ coverings with local
monodromies in $O$ and Galois groups contained in $G$, and we call
this semigroup an {\it admissible geometric degree $d$ covering
semigroup with local monodromies in $O$ and Galois group in $G$}.

The following statement is straightforward.

\begin{claim} \label{a=g} Any  admissible geometric degree $d$
covering semigroup ${\rm G\mathbb S_{d,\mathcal {\tilde
H}-equiv}}(G,O)$ is isomorphic over $G$ to the admissible algebraic
covering semigroup ${\rm \mathbb S_{\mathcal H-equiv}}(G,O)$, where
$\mathcal {H}=\{H_{n,p}\}_{n\geq 0,p\geq0}$ is the collection of
subgroups representing $\mathcal {\tilde H}=\{ \widetilde
H_{n,p}\}_{n\geq 0,p\geq0}$
 in
$Aut(\mathbb F^{n+2p})$.
\end{claim}

\subsection{Construction of Hurwitz spaces of marked coverings} \label{Hursp}
Here, we adapt Fulton's construction of Hurwitz spaces, see
\cite{F}, to the case of marked coverings.

Let $D\subset \overline F$ be an open disc in a projective
irreducible non-singular algebraic curve $\overline F$. We put
$F=\overline F\setminus D$, choose a point $q\in \partial F$, and fix an $n$-point set $B\subset F\setminus \partial F$.

Let us recall that for any surface $F$ its braid group on $n$
strands, $Br_n(F, \partial F)$, can be seen as fundamental group
$\pi_1(F^{(n)}\setminus \Delta)$, where $F^{(n)}$ is the symmetric
product of $n$ copies of $F$ and $\Delta$ is the discriminant locus.

Due to our assumptions, the fundamental group $\pi_1(F\setminus
B,q)$ is isomorphic to the free group $\mathbb F^{n+2p}$
where $p$ is the genus of $F$, and in such a way the braid group
$Br_n(F,\partial F)$, which acts naturally (the right action) on
$\pi_1(F\setminus B,q)\simeq\mathbb F^{n+2p}$, becomes
anti-isomorphic to the algebraic braid group $Br_{n,p}$ introduced
in Subsection \ref{adm}
(this is usually called {\it Artin presentation
theorem} and follows from a comparison of the actions
of the generators of these groups on $\mathbb F^{n+2p}\simeq \pi_1(F\setminus B,q)$).

To detail these identifications, let us pick a caudate skeleton
$S^{cdt}$ of $F$ the set of ends of whose tails is $B$. In notation
of subsection \ref{skeletons}, the choice of $S^{cdt}$ defines the
set $\{ \lambda_1,\mu_1,\dots, \lambda_p,\mu_p\}$ of free generators
of the group $\pi_1(F,q)$ and loops $\gamma_i$ around the points of
$B$, so that $\gamma_1,\dots, \gamma_n,\lambda_1,\mu_1,\dots,
\lambda_p,\mu_p$ are free generators of $\pi_1(F\setminus B,q)\simeq
\mathbb F^{n+2p}$ and $\gamma_1\dots \gamma_n[\lambda_1,\mu_1]\dots
[\lambda_p,\mu_p]=\partial F$ in $\pi_1(F\setminus B,q)$ (as usual
$\partial F$ is taken counter-clockwise).

The set $\{ \gamma_1,\dots, \gamma_n,\lambda_1,\mu_1,\dots,
\lambda_p,\mu_p\}$ is a frame of the free group $\pi_1(F\setminus
B,q)\simeq \mathbb F^{n+2p}$, which in accordance with notation of
Subsection \ref{adm} we denote by $\mathcal B_{{\bf{1}}}$. As to the
standard generators
$\sigma_1,\dots,\sigma_{n-1},\xi_{1,\lambda},\dots
,\xi_{p,\lambda},\xi_{1,\mu},\dots,
\xi_{p,\mu},\zeta_{1},\dots,\zeta_{p}$ of the algebraic braid group
$Br_{n,p}$, they turn in terms of the geometric braid group
$Br_n(F,\partial F)$ into H-moves, $\lambda$- and $\mu$-moves, and
braiding of a point around a handle (see Subsection
\ref{geom-strong}).

This anti-isomorphism $Br_n(F,\partial F)\to Br_{n,p}$ defines a
right action of $Br_n(F, \partial F)$ on the set of words
$W_{\mathcal B_{{\bf{1}}}}(G,O)=\bigcup_{f}W_{f,\mathcal
B_{{\bf{1}}}}$, where the union is taken over all equipped
homomorphisms $f:\mathbb F^{n+2p}\to (G,O)$ to an equipped group
$(G,O)$. Obviously, the subset $W^G_{\mathcal
B_{{\bf{1}}},{{\bf{1}}}}(G,O)$ of words in $W_{\mathcal
B_{{\bf{1}}}}(G,O)$ representing the elements of the semigroup
$\mathbb S(G,O)^G_{{\bf{1}}}$  is invariant under the action of
$Br_n(F,\partial F)$. Therefore this action  defines homomorphisms
$\omega=\omega_{n,p,(G,O)} : \pi_1(F^{(n)}\setminus \Delta)=
Br_n(F,\partial F)\to \mathcal S_{|W_{\mathcal B_{{\bf{1}}}}(G,O)|}$
and $\omega^G_{{\bf{1}}} : \pi_1(F^{(n)}\setminus \Delta)=
Br_n(F,\partial F)\to \mathcal S_{|W^G_{\mathcal B_{{\bf{1}}},
{{\bf{1}}}}(G,O)|}$ to the symmetric groups $\mathcal
S_{|W_{\mathcal B_{{\bf{1}}}}(G,O)|}$ and $\mathcal
S_{|W^G_{\mathcal B_{{\bf{1}}},{{\bf{1}}}}(G,O)|}$ acting,
respectively, on the sets $W_{\mathcal B_{{\bf{1}}}}(G,O)$ and
$W^G_{\mathcal B_{{\bf{1}}},{{\bf{1}}}}(G,O)$.

Put $F_0=F\setminus \partial F$. As it follows from Proposition
\ref{encoding}, the homomorphisms $\omega$ and $\omega^G_{{\bf{1}}}$
define a $|W_{\mathcal B_{{\bf{1}}}}(G,O)|$-sheeted  unramified
covering
$$
\theta_n=\theta_n(G,O): \widetilde{\text{HUR}}_{(G,O),n}(F)\to
F_0^{(n)}\setminus \Delta
$$
and a $|W^G_{\mathcal B_{{\bf{1}}},{{\bf{1}}}}(G,O)|$-sheeted
unramified covering
$$
\theta_n^G=\theta_n(G,O)^G_{{\bf{1}}}: \text{HUR}_{(G,O),n}(F)\to
F_0^{(n)}\setminus \Delta,
$$
respectively. Furthermore, there is a canonical embedding
$j:\text{HUR}_{(G,O),n}(F)\hookrightarrow
\widetilde{\text{HUR}}_{(G,O),n}(F)$ such that $\theta_n^G =\theta_n
\circ j$.  Moreover, the both coverings are marked (by words of
$W_{\mathcal B_{{\bf{1}}}}(G,O)$ and $W^G_{\mathcal
B_{{\bf{1}}},{{\bf{1}}}}(G,O)$, respectively) over the point in
$F_0^{(n)}\setminus \Delta$ represented by $B$, and $j$ respects the
markings.

According to the usual construction of covering spaces by means of
the groupoid of homotopy classes of paths, the covering space
$\widetilde{\text{HUR}}_{(G,O),n}(F)$ as a set is the set of pairs
$(B', f')$, where $B'\in F_0^{(n)}\setminus \Delta$ and $f' :
\pi_1(F\setminus B',q)\to G$ are epimorphisms such that the
conjugacy classes of their values at the loops around the points of
$B'$ belong to $O$ and $f'(\partial F)=\bf{1}$. We call
$\widetilde{\text{HUR}}_{(G,O),n}(F)$ the {\it Hurwitz space of
marked $n$-branched coverings of $F$ with equipped Galois group
$(G,O)$.} This construction being functorial, a choice of an
embedding $i:G\hookrightarrow \mathcal S_d$ provides an embedding of
$\widetilde{\text{HUR}}_{(G,O),n}(F)$ into
$\widetilde{\text{HUR}}_{d,n}=\widetilde{\text{HUR}}_{(\mathcal
S_d,\mathcal S_d\setminus \{ {{\bf{1}}}\}),n}(F)$, which we call the
{\it Hurwitz space of marked $n$-branched degree $d$ coverings of
$F$.}

The advantage of considering marked coverings is that the Hurwitz
spaces $\widetilde{\text{HUR}}_{d,n}$  come then with a universal
family of coverings, $\mathcal F_{d,n}\to
\widetilde{\text{HUR}}_{d,n}$.  Such a family can be obtained as
manifold completion ({see Subsection \ref{encoding-ramified}) of the
unramified covering of $U=\{ (p, B') : p\notin B'\}\subset F\times
(F_0^{(n)}\setminus\Delta) $, that is the covering defined by the
homomorphism $\pi_1(U, (q,B))\to \mathcal S_{I_d\times W_{\mathcal
B_{{\bf{1}}}}(G,O)}$, $G=\mathcal S_d$, $O=\mathcal S_d\setminus
\{1\}$ that sends the images of elements
$\varsigma\in\pi_1(F\setminus B, q)$ to permutations $(x,w)\mapsto
(\alpha(\varsigma)x,w)$ and the images of elements $\chi\in
\pi_1(F^{(n)}\setminus \Delta)$ to permutations $(k,w)\mapsto
(k,\omega(\chi)w)$.

\begin{claim} \label{mc}
The connected components of the spaces
$\widetilde{\text{HUR}}_{(G,O),n}(F)$ and $\text{HUR}_{(G,O),n}(F)$
are in one-to-one correspondence with the elements of the strong
covering semigroup $\mathbb S(G,O)$ and its subsemigroup $\mathbb
S(G,O)^G_{{\bf{1}}}$, respectively.
\end{claim}
\proof  By Proposition \ref{encoding} and due to existence of
universal families, the connected components of
$\widetilde{\text{HUR}}_{(G,O),n}(F)$ (respectively,
$\text{HUR}_{(G,O),n}(F)$) are in one-to-one correspondence with the
orbits of the action of $Br_n(F)$ on the set $W_{\mathcal
B_{{\bf{1}}}}(G,O)$ (respectively, $W^G_{\mathcal
B_{{\bf{1}}},{{\bf{1}}}}(G,O)$). Via the anti-isomorphism
$Br_n(F)\to Br_{n,p}$ and  due to the definition of the strong
covering semigroups, these orbits coincide with  the elements of
$\mathbb S(G,O)$.
 \qed \\

Let $F=\overline F\setminus D$ and $F'=\overline F\setminus D'$,
where $D'\subset D$ are open discs in a closed genus $p$ oriented
surface $\overline F$ without boundary such that the marked point
$q\in
\partial D\cap\partial D'$. Then there is a natural embedding
$j_{F,F'}:\text{HUR}_{(G,O),n}(F)\hookrightarrow
\text{HUR}_{(G,O),n}(F')$ which is compatible with the covering maps
and the embedding $j_{F,F'}:F^{(n)}_0\setminus \Delta\hookrightarrow
F'^{(n)}_0\setminus\Delta$. Let $i:G\hookrightarrow \mathcal S_d$ be
an embedding of a group $G$ such that its image acts transitively on
$I_d$. If a word $w_{f_*}$ represents an element of $\mathbb
S(G,O)^{G}_{{\bf{1}}}$, then, first, the covering space $E$ of the
$d$ sheeted marked covering $f:E\to F$ is connected, and, second,
the covering $f$ can be extended uniquely to a $d$ sheeted marked
(at $q\in \overline F$) covering $f:\overline E\to \overline F$
unbranched at the points of $D$. The embeddings
$\text{HUR}_{(G,O),n}(\overline F\setminus D_i)\hookrightarrow
\text{HUR}_{(G,O),n}(\overline F\setminus D_{i+1})$ corresponding to
an infinite sequence of open discs
$$\dots \subset D_{i+1}\subset D_i\subset \dots \subset D_{1}$$ such that
$\cap_{i=1}^{\infty} D_i=\emptyset$ and $q\in \partial D_i$ for all
$i$ define an unramified covering $\theta_n(G,O)^G_{{\bf{1}}}:
\text{HUR}_{(G,O),n}(\overline F,q)\to (\overline F\setminus
q)^{(n)}\setminus \Delta$ the covering space of which is called the
{\it Hurwitz space} of marked (at a point $q\in \overline F$)
coverings of a projective algebraic
 curve $\overline F$ with equipped Galois group $(G,O)$ and branched at $n$
points.

\subsection{Proof of Theorems \ref{cor1}, \ref{TH2}, and \ref{TH1}}
Due to Proposition \ref{geom=alg}: Theorem \ref{cor1} follows from Corollary \ref{cor2};
Theorem \ref{TH2} follows from  Claim \ref{mc}, Theorems \ref{Cuniq3},  and \ref{Cuniq4};
and Theorem \ref{TH1}  follows from  Claim \ref{mc} and Corollary \ref{cor1}.

 \ifx\undefined\bysame
\newcommand{\bysame}{\leavevmode\hbox to3em{\hrulefill}\,}
\fi


\begin{thebibliography}{McD-Sa-2}

\def\entry#1#2#3#4\par{\bibitem[#1]{#1}
{\textsc{#2 }}{\sl{#3}} #4\par\vskip2pt}

\def\noentry#1#2#3#4\par{}

\bibitem
{BE} {I.~Berstein and A.L.~Edmonds:} {\it On the classification of
generic branched coverings of surfaces.} Illinois J. Math. Volume
28, Issue 1 (1984), 64 -- 82.

\bibitem 
{F}%
{W. Fulton:} {\it Hurwitz schemes and irreducibility of moduli of
algebraic curves.} Ann. of Math., {\bf 90:3} (1969), 542 -- 575.

\bibitem 
{Cl} {A. Clebsch:} {\it Z$\ddot{u}$r Theorie der Riemann'schen
Fl$\ddot{a}$che.} Math. Ann., 6 (1872), 216 -- 230.

\bibitem
{FaN} {E.~Fadell , L.~P. Neuwirth:}{\it Configuration spaces,} Math.
Scand. 10 (1962), 111 -- 118.

\bibitem
{Fox} {R.~H.~Fox:} {\it Covering Spaces with Singularities.} In
Algebraic Geometry and Topology. A Symposium in honor of S.
Lefschetz. Editors: Fox et al. Princeton Univ. Press (1957) pp.243
-- 257.

\bibitem
{FoN} {R.~H.~Fox, L.~P. Neuwirth:}{\it The braid groups,} Math.
Scand. 10 (1962), 119 -- 126.

\bibitem
{FV} {M.D. Fried and H. V$\ddot{o}$lklein:} {\it The inverse Galois
problem and rational points on moduli space.} Math. Ann., 290,
(1991), 771 -- 800.

\bibitem
{GHS}{T.~Graber, J.~Harris, J.~Starr:}
{\it   A note on Hurwitz schemes of covers of a positive genus curve.}
arXiv:math/0205056

\bibitem
{H} {A. Hurwitz:} {\it Ueber Riemann'she Fl$\ddot{a}$chen mit
gegebenen Verweigugspunkten.} Math. Ann., 39, (1981), 1 -- 61.

\bibitem
{Ka} {V. Kanev:} {\it Irreducibility of Hurwitz spaces.} arXiv:
math/0509154v1 [math.AG] 7 Sep 2005.

\bibitem 
{Ku} {Vik.S. Kulikov:} {\it Hurwitz curves.} UMN {\bf 62:6} (2007),
3 -- 86.

\bibitem 
{Ku1} {Vik.S. Kulikov:} {\it Factorization semigroups and
irreducible components of Hurwitz space,}  Izv. Math. {\bf 75:4}
(2011), 711 -- 748. 

\bibitem 
{Ku2} {Vik.S. Kulikov:} {\it Factorization semigroups and
irreducible components of Hurwitz space. II,} (accepted in Izv.
Math.; primiry version is in arXiv:1011.3619).

\bibitem 
{Ku3} {Vik.S. Kulikov:} {\it Factorizations in finite groups,}
arXiv:1105.1939 (submitted in Sb. Math.).

\bibitem 
{KK} {V. Kharlamov and Vik.S. Kulikov: } {\it On braid monodromy
factorizations.} Izv. Math. 67:3 (2003), 499 -- 534.

\bibitem 
{P} {A. N. Protopopov: } {\it Topological classification of branched coverings of the two-dimensional sphere.}
Journal of Mathematical Sciences. 52:1 (2003), 2832--2846.

\bibitem {V} {F. Vetro:} {\it Irreducibility of Hurwitz spaces for coverings with one special fibre.} Indag. Math. (N.S.),
vol. 17 (2006), no. 1, 115 -- 127.


\bibitem {W} {B. Wajnryb:} {\it Orbits of Hurwitz action for
coverings of a sphere with two special fibres.} Indag. Math. (N.S.),
vol. 7 (1996), no. 4, 549 -- 558.


\end{thebibliography}
\end{document}